\RequirePackage{scrlfile}
\ReplacePackage{scrpage2}{scrlayer-scrpage}
\documentclass[12pt,reqno]{amsart}
\newcommand{\mycomment}[1]{} 
\usepackage{style/amsart-uppercase-title}
\usepackage[utf8]{inputenc}
\usepackage[backend=bibtex,giveninits=true,url=true,style=alphabetic,citestyle=alphabetic,isbn=false,backref]{biblatex}
\usepackage{amssymb,amsfonts,amsmath,amsthm}
 \usepackage{geometry}
 \usepackage{empheq}
\usepackage[most]{tcolorbox}
\newtcbox{\mymath}[1][]{
  nobeforeafter, math upper, tcbox raise base,
    enhanced, colframe=blue!30!black,
    colback=green!30, boxrule=1pt,
    #1}
\allowdisplaybreaks

\def\upperrel#1#2{\mathrel{\mathop{#2}\limits^{#1}}}

\usepackage{style/mypub}
\usepackage{style/mysymb}
\usepackage{amsbsy}

\usepackage{refcount}
\usepackage{style/fullref}
\usepackage{style/autoref-addons}
\usepackage{thmtools}
\usepackage{style/mythmtools-setup}
\usepackage{enumitem}
\setlist[enumerate]{label=\arabic*.)}

\usepackage{xr-hyper}
\usepackage[colorlinks,citecolor=blue,filecolor=red,pdfnewwindow=true]{hyperref}
\renewcommand\equationautorefname{Inequality}

\usepackage{subfiles}
\begin{document}
 \title[A global existence theorem for the Euler-Nordström
 system]{Global existence of solutions to the irrotational
   Euler-Nordström equations with a positive cosmological constant: The gravitational field equation}

 \author[U.~Brauer]{Uwe Brauer}
 \address{
   Uwe Brauer Departamento de Matemática Aplicada\\ Universidad
     Complutense Madrid 28040 Madrid, Spain} \email{oub@mat.ucm.es}
 \thanks{U.~B.~gratefully acknowledges support from Grant PID2022-137074NB MINECO, Spain and UCM-GR17-920894.}
 \author[L.~Karp]{Lavi Karp}
   \address{
     Lavi Karp
     Department of Mathematics\\ Braude College of Engineering \\
     Snunit 51 St.,  
     Karmiel 2161002\\
     Israel}
   \email{karp@braude.ac.il}
\begin{abstract}

Our aim is to establish the global existence of classical solutions to
the nonlinear irrotational Euler–Nordström system, which incorporates
a linear equation of state and a cosmological constant. 
In this setting, gravitation is described by a single scalar field
satisfying a specific semilinear wave equation. 
We restrict attention to spatially periodic perturbations of the
background metric and therefore study this equation on the
three-dimensional torus $\setT^3$, working within the Sobolev spaces
$H^m(\setT^3)$.

We begin by analysing the Nordström equation in isolation, with a
source term generated by an irrotational fluid obeying a linear
equation of state. 
This separation is motivated by the fact that such a fluid produces a
source term containing a nonlinear contribution of fractional order.

To obtain a global solution for the gravitational field, the
fractional-order nonlinearity $(1+u)^\mu$, with $\mu\in\mathbb{R}$,
must remain smooth throughout the evolution. 
This condition, in turn,  requires  that $u$ remain small for all time. 
We ensure this by introducing a suitably chosen energy functional.
We also prove that, asymptotically, the solutions tend to a constant.

\end{abstract}
\maketitle{}
\tableofcontents
\renewenvironment{proof}[1][\proofname]{{\bfseries Proof of #1.}}{\hfill$\blacksquare$}

\makeatletter 
  \def\plain@equationautorefname{Equation}
\renewenvironment{proof}[1][\proofname]{{\bfseries Proof of #1.}}{\hfill$\blacksquare$}
\section{Introduction}
\label{sec:introduction}

The aim of this paper is to prove the global existence and uniqueness
of classical solutions, in the Sobolev spaces $H^m(\setT^3)$, for a
nonlinear wave equation with damping. 
The equation arises from the nonlinear Nordström--Euler system, which
models an irrotational relativistic fluid with a linear equation of
state.

We study small perturbations of a vacuum background solution. 
In this setting, the irrotational Euler flow reduces to a coupled
system of wave equations: a scalar equation for the gravitational
field, whose source term contains a fractional power, and an
acoustical equation governing the streamlines. 
The fractional-power nonlinearity poses substantial technical
obstacles to proving global existence. 
To address them, we introduce a carefully designed energy functional,
defined in \autoref{M-eq:decaying-energy-estimates-nordstroem:3}. 
This functional is the central tool in our proof of global existence
for the Nordström gravitational field equation.

\subsection{The field equations with cosmological constant and the
  background solutions}
\label{sec:fields-equations}

The first fully relativistic, consistent, theory of gravitation was a
scalar theory developed by Nordström
\cite{nordstrom13:_zur_theor_gravit_stand_relat}, where the
gravitational field is described by a nonlinear hyperbolic equation
for the scalar field $\phi$.
Although the theory is not in agreement with observations it provides,
due to its nonlinearity, some interesting mathematical challenges. 
Surprisingly, this theory has never been mathematically investigated,
although its linear version coupled to the Euler equations has been
studied by Speck \cite{Speck_09} and coupled to the Vlasov equation by
Calogero \cite{calogero03:_spher} and
others~\cite{Felix_Antonio_Calogero_2014},
\cite{Fajman_Jeremie_Jacques-2021}, \cite{Wang_2021},
\cite{Calogero_Rein_2003} and \cite{Calogero_Rein_2004}.

We follow here the geometric reformulation provided by Einstein-Fokker
\cite{einstein14:_nords_gravit_stand_differ%
} and will use the Euler equations as a matter model. 
See also Straumann~\cite[Chap.~2.]{STRA2%
} for a modern representation of that theory.

The basic idea of this theory is that the physical metric
$ g_{\alpha\beta}$ is related to the Minkowski metric
$\eta_{\alpha\beta}$ by the following conformal transformation.
\begin{equation}
  \label{eq:Euler-Nordstrom:1}
   g_{\alpha\beta}=\phi^2\eta_{\alpha\beta}, 
\end{equation}
where $\eta_{\alpha\beta}=\text{diag}(-1,1,1,1)$.

The Nordström field equations with a cosmological constant takes the form 
\begin{equation}
  \label{eq:Nordstrom:3}
  \square\phi =- \frac{1}{6} T\phi^{3} -\Lambda \phi, 
\end{equation}
where 
\begin{equation}
\label{eq:Alembert}
 \square \eqdef 
\eta^{\alpha\beta}\partial_{\alpha}\partial_{\beta}= -\partial_{t}^{2}+\Delta.
\end{equation}

Here $T$ denotes the  trace of the fundamental energy tensor 
$T=g_{\alpha\beta}T^{\alpha\beta}$, which  for a perfect fluid is given by 
\begin{equation}
  \label{eq:section1-intro:6}
   T^{\alpha\beta}= \left(  \epsilon + p \right) u^{\alpha} u^{\beta}+  p g^{\alpha\beta},
\end{equation}
where $\epsilon$ denotes the energy density, $p$ the pressure 
and $u^{\alpha}$ is the unit timelike vector which satisfies
\begin{equation}
  \label{eq:section1-intro:7}
   g_{\alpha\beta} u^{\alpha} u^{\beta}=-1.
 \end{equation}
Moreover, we assume   a linear equation of state
\begin{equation}
  \label{eq:nordstroem-field:9}
  p=K\epsilon, \qquad K>0
\end{equation}
and the fluid to be irrotational which we will discuss in \fullref{sec:plugin-fluid}.

\subsubsection{Cosmological setting}
\label{sec:cosmological-setting}

We assume an isotropic and homogeneous vacuum background solution,
which we denote by $\mathring{\phi}$.
Homogeneity implies that the function $\mathring \phi$ depends just on
$t$, while the fact that the solution describes vacuum leads to the
conclusion that $T\equiv 0$.
Therefore \autoref{eq:Nordstrom:3} reduces to
\begin{equation}
  \label{eq:section1-intro:1B}
  -\frac{d^2}{dt^2}\mathring{\phi}=-\Lambda\mathring{\phi}.
\end{equation}
For convenience we set
\begin{equation}
  \label{eq:nordstroem-field:39}
  \Omega^2=\Lambda>0.
\end{equation}
Therefore this differential equation has a general solution of the form
\begin{equation*}
  \mathring{\phi}=Ae^{\Omega t}+Be^{-\Omega t}.
\end{equation*}

Since we want our solution to mirror the so-called flat de Sitter solution in general relativity (see, for example,
\cite[Chap.~V]{Choquet-Bruhat_09%
}
so that both $\mathring{\phi}$ and $\frac{d}{dt}\mathring{\phi}$ are positive, we choose
\begin{equation}
\label{eq:section1-intro:2B}
\mathring{\phi}(t)= \mathord{e^{\Omega t}}
\end{equation}
as the background solution.

\subsubsection{Deviations}
\label{sec:deviations}

We now study small deviations from the background solution
$\mathring \phi$.
So we make the following Ansatz
\begin{align}
  \label{eq:phi}
  \phi&=\mathring{\phi}+\Psi=\mathord{e^{\Omega t}} +\Psi,\\
\label{eq:nordstroem-field:8}
  T&= {\mathring{T}} + \stackrel{1}{T}
\end{align}
where $\Psi$ denotes the deviation from the background.
We assume the background solution to be vacuum, so that
$\mathring{T}=0$. It then follows, by taking  the trace
of the energy-momentum tensor, with respect to
the metric \autoref{eq:Euler-Nordstrom:1} and the linear equation of state
\autoref{eq:nordstroem-field:9}, that $\stackrel{1}{T}$ takes the form
\begin{equation}
  \label{eq:nordstroem-field:7}
  \stackrel{\scriptscriptstyle 1}{{T}}
  = 3\stackrel{\scriptscriptstyle 1}{p}-\stackrel{\scriptscriptstyle 1}{\epsilon}
  = 3\stackrel{\scriptscriptstyle 1}{p}-\frac{1}{K}\stackrel{\scriptscriptstyle 1}{p}= 
\left( 3- \frac{1}{K} \right)\stackrel{\scriptscriptstyle 1}{p}.
\end{equation}
We might drop the $\stackrel{\scriptscriptstyle 1}{{T}}$  in the following 
whenever it is convenient.
The deviation $\Psi$ satisfies the following equation
\begin{equation}
  \square \phi=\square(\mathord{e^{\Omega t}}+\Psi)=-\Omega^2 \mathord{e^{\Omega t}} + \square \Psi
  =-\frac{1}{6} \stackrel{1}T(\mathord{e^{\Omega t}} +\Psi)^3-\Omega^2(\mathord{e^{\Omega t}} +\Psi).
\end{equation}
Thus $\Psi$ satisfies the initial value problem
\begin{equation}
  \label{eq:psi}
  \left\{
    \begin{array}{l}
      \square \Psi =-\frac{1}{6}\stackrel{1}{T}(\mathord{e^{\Omega t}}+\Psi)^3-\Omega^2\Psi\\
      \Psi(0,x)=\Psi_0(x), \\
      \partial_t\Psi(0,x)=\Psi_1(x)
    \end{array}\right..
\end{equation}
Our goal is to show global existence of classical solutions for
\autoref{eq:psi} demanding a small source term $T$ and small
initial data.

Note that if $\Psi$ is small, then
\begin{equation*}
  (\mathord{e^{\Omega t}} +\Psi)^3\sim e^{3\Omega t},
\end{equation*}
and this term growths very rapidly and might prevent that the solution
exists for all time.
So in order to achieve the desired asymptotic behaviour of $\Psi$, we
multiply $\phi$ by $\mathord{e^{-\Omega t}}$, then from \autoref{eq:phi} we
conclude that
\begin{equation*}
  \mathord{e^{-\Omega t}} \phi=1+\mathord{e^{-\Omega t}}\Psi ,
\end{equation*}
and therefore we set
\begin{equation}
\label{eq:section1-intro:5}
  \Theta\overset{\mbox{\tiny{def}}}{=} \mathord{e^{-\Omega t}} \Psi,
\end{equation}
or
\begin{equation}
  \label{eq:nordstroem-field:17}
  \phi = e^{\Omega t}(1+\Theta).
\end{equation}

In order to obtain a suitable equation  for $\Theta$ we proceed as follows:
\begin{align*}
  \partial_t \Theta &= \partial_t (e^{-\Omega t} \Psi)=e^{-\Omega t} \partial_t\Psi-\Omega e^{-\Omega t} \Psi=e^{-\Omega t} \partial_t\Psi-\Omega \Theta,\\
  \partial_t^2 \Theta &= e^{-\Omega t}\partial_t^2\Psi-2\Omega \partial_t \Psi+\Omega^2 e^{-\Omega t}\Psi= e^{-\Omega t}\partial_t^2\Psi-2\Omega \left(\partial_t \Theta+\Omega\Theta)\right) +\Omega^2\Theta \\
       & = e^{-\Omega t}\partial_t^2\Psi-2\Omega \partial_t \Theta-\Omega^2\Theta,
\end{align*}
or
\begin{equation}
\label{eq:nordstroem-field:3}
  -e^{-\Omega t}\partial_t^2\Psi=-\partial_t^2\Theta-2\Omega \partial_t   \Theta -\Omega^2\Theta.
\end{equation}
Inserting it in \autoref{eq:psi}, we get 
\begin{align*}
  -\partial_t^2\Theta-2\Omega \partial_t   \Theta -\Omega^2\Theta+e^{-\Omega t}\Delta \Psi 
  &  = - \frac{1}{6}e^{-\Omega t} \stackrel{\scriptscriptstyle 1}{{T}}(e^{\Omega t}+\Psi)^{3}-\Omega^2e^{-\Omega t}\Psi \\
  & = - \frac{1}{6}e^{2\Omega t} \stackrel{\scriptscriptstyle 1}{{T}}(1+\Theta)^{3}-\Omega^2\Theta.
\end{align*}

Thus, the resulting equation for $\Theta$ takes the form
\begin{equation}
\label{eq:nordstroem-field:43}
 \square \Theta   -2 \Omega \partial_{t} \Theta=- \frac{1}{6}e^{2\Omega t} \stackrel{\scriptscriptstyle 1}{{T}}(1+\Theta)^{3},
\end{equation} 
where $\square$ is the d'Alembert operator defined in \autoref{eq:Alembert}.

\subsubsection{A irrotational fluid with a linear equation of state as a source}
\label{sec:plugin-fluid}

Let us recall the basic ideas following
\cite[Chap. IX, Sec. 10]{Choquet-Bruhat_09%
}.
For a  barotropic fluid, the Euler equations are derived from the
conservation law
\begin{equation}
  \label{eq:section1-intro:3}
  \nabla_{\alpha}T^{\alpha\beta}=0
\end{equation}
and take the form
\begin{align}
  \label{eq:chquet-bruhat-irrotational-fluid:2}
  (\epsilon+p) \nabla_{\alpha} u^{\alpha}+u^{\alpha} \partial_{\alpha} \epsilon=0
\end{align}
and
\begin{align}
    \label{eq:chquet-bruhat-irrotational-fluid:3}
  (\epsilon+p) u^{\alpha} \nabla_{\alpha} u^{\beta}+\left(g^{\alpha \beta}+u^{\alpha} 
u^{\beta}\right) \partial_{\alpha} p=0.
\end{align}
Now the index of the fluid is given by 
\begin{equation}
  \label{eq:chquet-bruhat-irrotational-fluid:6}
  f(p)\upperrel{\rm def}{=}\exp \int \frac{d p}{\epsilon(p)+p}
\end{equation}
and the dynamic of the 4-velocity vector $u^\alpha$  by
\begin{equation}
\label{eq:choquet-bruhat-irrotational-fluid:6}
  C_{\alpha}\upperrel{\rm def}{=}f u_{\alpha}, 
\end{equation}
hence
\begin{equation}
 g^{\alpha\beta}C_{\alpha} C_{\beta}=-f^{2}.
\end{equation} 
This function $f$ satisfies 
\begin{align}
  \label{eq:chquet-bruhat-irrotational-fluid:8}
  \partial_{\alpha} f& \equiv \frac{\partial f}{\partial x^{\alpha}}=f \quad
\frac{\partial_{\alpha} p}{\epsilon+p}
\end{align}
and
\begin{align}
  \label{eq:chquet-bruhat-irrotational-fluid:9}
  \partial_{\alpha} &f=-f^{-1} C^{\beta} \nabla_{\alpha} C_{\beta}.
\end{align}
Let
  \begin{math}
    \mu_{p}^{\prime} \equiv \frac{\partial \mu}{\partial p} 
  \end{math}.
   According to
\cite[Thm.~10.1]{Choquet-Bruhat_09%
}
the Euler equations \eqref{eq:chquet-bruhat-irrotational-fluid:2}
and \eqref{eq:chquet-bruhat-irrotational-fluid:3} are equivalent to the following 
equations:
  \begin{align}
    \label{eq:chquet-bruhat-irrotational-fluid:10}
    C^{\alpha}\left(\nabla_{\alpha} C_{\beta}-\nabla_{\beta} C_{\alpha}\right)=0
\end{align}
and
\begin{align}
    \label{eq:chquet-bruhat-irrotational-fluid:11}
    \nabla_{\alpha} C^{\alpha}+\left(\mu_{p}^{\prime}-1\right) \frac{C^{\alpha} 
C^{\beta}}{C^{\lambda} C_{\lambda}} \nabla_{\alpha} C_{\beta}=0 .
  \end{align}

  \begin{defn}[Irrotational fluid flow]
    \label{def:section1-intro:1}
  The vorticity of a flow is an anti symmetric tensor 
  \begin{equation}
\label{eq:section1-intro:12}
   \Omega_{\alpha\beta}\equiv \nabla_{\alpha} C_{\beta}-\nabla_{\beta} C_{\alpha}.
  \end{equation} 
    A fluid flow with zero vorticity is called irrotational.
  \end{defn}
 
  The trajectories of irrotational flow are locally orthogonal to
  hypersurfaces. 
  The equation
  \begin{equation*}
    \nabla_{\alpha} C_{\beta}-\nabla_{\beta} C_{\alpha}=0
  \end{equation*}
  implies that there exists on a spacetime,
  at least locally, a function $\Phi$ such that
    \begin{equation}
      \label{eq:choquet-bruhat-irrotational-fluid:5}
      C_{\alpha}=\partial_{\alpha} \Phi.
    \end{equation}

    Hence for the irrotational Euler flow equations
    \eqref{eq:chquet-bruhat-irrotational-fluid:10} and
    \eqref{eq:chquet-bruhat-irrotational-fluid:11} are reduces to the
    quasi-linear wave equation
\begin{equation}
  \label{eq:section1-intro:4}
  \nabla_{\alpha} \partial^{\alpha} \Phi+\left(\mu_{p}^{\prime}-1\right) \frac{\partial^{\alpha} \Phi\partial^{\beta}\Phi}{\partial^{\lambda}\Phi \partial_{\lambda}\Phi} \nabla_{\alpha} \partial_{\beta}\Phi=0.
\end{equation}

This equation can be simplified further, but we skip this for the moment.    

With these concepts at hand we claim the following.

\begin{prop}[The source term of an irrotational fluid with an linear
  equation of state]
  \label{prop:section1-intro:1}
  For an irrotational fluid, defined by equations
  \eqref{eq:choquet-bruhat-irrotational-fluid:6},
  \eqref{eq:choquet-bruhat-irrotational-fluid:5}, with an linear
  equation of state \eqref{eq:nordstroem-field:9}, and a fluid index
  given by \autoref{eq:chquet-bruhat-irrotational-fluid:6}, the right
  hand side of \autoref{eq:nordstroem-field:43} takes the form
  \begin{equation}
    \label{eq:section1-intro:8}
    - \frac{1}{6}e^{2\Omega t} \stackrel{\scriptscriptstyle 1}{{T}}(1+\Theta)^{3}
    = -\frac{1}{6}   \left( 3-\frac{1}{K}\right)
    {e^{\frac{(K-1)\Omega}{K} t}}
    \left[ \left( \partial_{t}\Phi  \right)^{2} - \left\vert \nabla\Phi \right\vert^{2}   \right]^{\frac{1+K}{2K}}
    \quad (1+\Theta)^{\left(\frac{2K-1}{K} \right)},
  \end{equation}
  where
  \begin{equation}
    \label{eq:section1-intro:14}
    \nabla\Phi:=(\partial_{1}\Phi, \partial_{2}\Phi, \partial_{3}\Phi) \qquad  \left\vert \nabla\Phi \right\vert^{2}:=\delta^{ab} \partial_{a}\Phi\partial_b\Phi
  \end{equation}
\end{prop}

\begin{proof}[~\protect\autoref{prop:section1-intro:1}]

  For a linear equation of state \eqref{eq:nordstroem-field:9} the
  fluid index defined by \autoref{eq:chquet-bruhat-irrotational-fluid:6}
  $f$ takes the form
  \begin{align}
    \label{eq:choquet-bruhat-irrotational-fluid:20}
    f(p)=\exp\left(\int \frac 
    {dp}{\left(\frac{1}{K}+1\right)p}\right)=p^{\frac{K}{1+K}},
  \end{align}
  hence $f$ can be expressed as follows
  \begin{align}
    \label{eq:nordstroem-field:50}
    p=f^{\frac{1+K}{K}}
  \end{align}
  and therefore the trace of the energy momentum tensor takes the form
  \begin{equation}
    \label{eq:nordstroem-field:51}
    \stackrel{\scriptscriptstyle 1}{{T}}= \left( 3- \frac{1}{K} \right)p=\left( 3- 
      \frac{1}{K} \right)f^{\frac{1+K}{K}}.
  \end{equation}
  On the other hand, by 
  \autoref{eq:choquet-bruhat-irrotational-fluid:6} and by considering
  the conformal metric given by \autoref{eq:Euler-Nordstrom:1},
  \begin{align}
    \label{eq:section1-intro:2}
    -f^{2}&=  g^{\alpha\beta} C_{\alpha} C_{\beta}  = \phi^{-2}\eta^{\alpha\beta}      
           C_{\alpha} C_{\beta}.                         
  \end{align}
  So by the existence of the potential $\Phi$, 
  \autoref{eq:choquet-bruhat-irrotational-fluid:5}, we can express $f$ by 
  \begin{align}
    \label{eq:section1-intro:1}
    f=  \phi^{-1} \left[ \left( \partial_{t}\Phi  \right)^{2} - \left\vert \nabla\Phi \right\vert^{2}   \right]^{\frac{1}{2}}
  \end{align}
  or
  \begin{align}
    \label{eq:section1-intro:9}
    f^{\frac{1+K}{K}}&=  \phi^{-\frac{1+K}{K}}
                       \left[ \left( \partial_{t}\Phi  \right)^{2} - \left\vert \nabla\Phi \right\vert^{2}   \right]^{\frac{1+K}{2K}}.
  \end{align}
  So inserting  \autoref{eq:section1-intro:9} into 
  \autoref{eq:nordstroem-field:51} we end up with
  \begin{equation}
    \label{eq:nordstroem-field:55}
    \stackrel{\scriptscriptstyle 1}{{T}}=\left( 3- \frac{1}{K} \right)  f^{\frac{1+K}{K}} 
    = \phi^{-(\frac{1+K}{K})} \left[ \left( \partial_{t}\Phi\right)^{2} - \left\vert 
        \nabla\Phi \right\vert^{2} \right]^{\frac{1+K}{2K}}.
  \end{equation}

  Using  \autoref{eq:nordstroem-field:17} we see that 
  \autoref{eq:nordstroem-field:55} results in
  \begin{equation}
    \label{eq:section1-intro:13}
    - \frac{1}{6}e^{2\Omega t} \stackrel{\scriptscriptstyle 1}{{T}}(1+\Theta)^{3}
    = -\frac{1}{6}   \left( 3-\frac{1}{K}\right)
    {e^{\frac{(K-1)\Omega}{K} t}} 
    \left[ \left( \partial_{t}\Phi  \right)^{2} - \left\vert \nabla\Phi \right\vert^{2}   \right]^{\frac{1+K}{2K}}
    \quad (1+\Theta)^{\left(\frac{2K-1}{K} \right)}
  \end{equation}
  which finishes the proof.
\end{proof}

\begin{rem}[On the character of the source term in \protect \autoref{eq:section1-intro:13}]
  \label{rem:nordstroem-field:1}
  The constant $K$, defined by the equation of state $p=K\epsilon$,
  must satisfy $0<K<1$. 
  This restriction fixes the range of exponents in the three terms on
  the right-hand side of \autoref{eq:section1-intro:13}:
  \begin{enumerate}
    \item The exponent of the second term ranges as
    $1<\frac{1+K}{2K}<\infty$, though the full expression is assumed
    to be small; see \fullref{rem:section1-intro:2}. 
    \item The last exponent satisfies $0<\frac{2K-1}{K}<1$.
    For its treatment, see \fullref{M-prop:fractional}.
  \end{enumerate}
\end{rem}

\subsection{The Cauchy problem for the Nordström equation for the gravitational field}
\label{sec:cauchy-problem}

For the sake of clarity, we adjust the notation as follows:
\begin{align}
  \label{eq:section1-intro:10}
  u      &\upperrel{\rm def}{=}   \Theta, \\
  \label{eq:section1-intro:11}
  a(t,x) &\upperrel{\rm def}{=}\frac{1}{6}   \left( 3-\frac{1}{K}\right) \left[ \left(
          \partial_{t}\Phi  \right)^{2} - \left\vert \nabla\Phi \right\vert^{2}
          \right]^{\frac{1+K}{2K}},   \\
  \varkappa      &\upperrel{\rm def}{=} \left(\frac{(1-K)\Omega}{K}\right)
\end{align}
and
\begin{align}
    \mu\upperrel{\rm def}{=}\frac{2K-1}{K}.
\end{align}
Please note, that $K$ has the range  $0<K<1$, therefore $\varkappa$ is always strictly positive.

With this notation in place, we study the following Cauchy problem:
\begin{subequations}
  \begin{align}
    \label{eq:wave:1}
    \partial_t^2 u -\Delta u  &=- 2\Omega \partial_t u + \mathord{e^{-\varkappa t}} a(t,x)(1+u)^{\mu}\\
    \label{eq:wave:2}
 &  (u(0,\cdot),\partial_t u(0,\cdot)) =(u_0,u_1),
  \end{align}
  where $0<\varkappa$ and $\mu\in\setR$.
\end{subequations}
We note that the deviations are required to be spatially periodic. 
For this reason, equations \eqref{eq:wave:1}--\eqref{eq:wave:2} are
formulated in the Sobolev spaces $H^m(\setT^3)$.

\begin{rem}[On the smallness assumption of $a(t,x)$]
  \label{rem:section1-intro:2}
  The physical regime under consideration is that of a vacuum solution
  subject to a small but finite perturbation. 
  For a linear equation of state, this entails a small pressure $p$,
  which in turn implies a small trace of the energy--momentum tensor. 
  From \autoref{eq:choquet-bruhat-irrotational-fluid:20},
  \autoref{eq:section1-intro:1}, and \autoref{eq:section1-intro:11},
  it follows that $a(t,x)$ must itself be small.
\end{rem}

To secure a global solution of the gravitational field equation, the
term $(1+u)^\mu$ must remain smooth throughout the evolution. 
This, in turn, requires that $u$ stay small at all times. 
To achieve this, we construct the energy functional by
\autoref{M-eq:decaying-energy-estimates-nordstroem:3}. 
Beyond controlling the time derivative and spatial gradient of the
solution, it also bounds $\|u\|_{H^m(\setT^3)}$. 
By keeping the energy uniformly small, we ensure that $u$ itself
remains small.

In an earlier paper \cite{Brauer_Karp_23}, we studied a related
semi-linear equation with a cubic nonlinearity, $(1+u)^3$, which is
smooth for all values of $u$. 
That analysis was possible because the equation was not coupled to the
Euler system \eqref{eq:section1-intro:3}. 
There, we relied on the machinery of symmetric hyperbolic systems, but
the resulting energy estimate did not control $u$. 
With fractional powers, by contrast, such control is indispensable,
both for global existence and for the long-time asymptotic behaviour
of solutions.

There is a substantial body of work on the global existence of
semi-linear wave equations with small initial data (see, for example,
\cite{Ebert_Reissig_18}). 
Related results are available for wave equations on a de Sitter
spacetime \cite{Yagdjian_Galstian_08%
}.

The problem considered here has a different structure. 
The right-hand side is of the form $e^{-\varkappa} a(t,x)(1+u)^\mu$, and the
smallness assumption falls not on the initial data but on the
coefficient $a(t,x)$. 
To our knowledge, this class of problems has not been addressed in the
existing literature.

Having set out the problem, we describe the structure of the paper and
summarise our main results.

In \fullref{M-sec:math-prel} we introduce the necessary mathematical
tools, such as Sobolev spaces on the torus $\setT^3$ and estimates for
nonlinear functions.

We then turn, in \fullref{M-sec:energy-estimates}, to our principal
analytical device: the energy estimates.

Our main result follows in \fullref{M-sec:main-results}, where it is
presented together with its proof and two essential lemmas.

The proofs of these lemmas are given in
\fullref{M-sec:proof-fullr-lem:s}.

Finally, in \fullref{M-sec:asymptotic-behaviour}, we establish a
result concerning the asymptotic behaviour of the global solution.

\makeatletter 
  \def\plain@equationautorefname{Inequality}

\renewenvironment{proof}[1][\proofname]{{\bfseries Proof of #1.}}{\hfill$\blacksquare$}

\section{Mathematical preliminaries}
\label{sec:math-prel}

\subsection{Sobolev spaces on the torus $\setT^3$}
\label{sec:non-homog-sobol}
We begin by fixing a convention for the indices. 
Greek indices such as $\alpha, \beta$ range over the values ${0,1,2,3}$. 
For non-negative integers $(n_1, n_2, n_3)$, we define a spatial
multi-index as $\vec{\alpha}= (n_1, n_2, n_3)$ and
$\partial_{\vec{\alpha}}=(\partial_1^{n_1}, \partial_2^{n_2}, \partial_3^{n_3})$.

For the convenience of the reader we present here the definition of the Sobolev space we use in this article.
\begin{defn}[Sobolev spaces over the torus]
  \label{def:section2-preliminaries:1}
  The Sobolev spaces over the torus $\setT^{3}$ (with the flat metric) for
  a nonnegative integer $m $ are defined as a completion of
  $C^{m}(\setT^3)$ functions under the norm
  \begin{equation}
    \label{eq:Sobolev-norm}
    \| u\|_{H^m(\mathbb T^3)}^2=\sum_{|\mathord{\vec{\alpha}}|\leq m} \|\partial_{\mathord{\vec{\alpha}}}  u\|_{L^2(\mathbb 
T^3)}^2 =\sum_{|\mathord{\vec{\alpha}}|\leq m}\int_{\mathbb T^3} |\partial_{\mathord{\vec{\alpha}}}  u(x)|^2d^{3}x,
  \end{equation}
  here $\vec{\alpha}$ denotes a spatial multi-index. 
  We denote them by $H^m(\mathbb T^3)$, and throughout the paper, we
  will use the abbreviation $H^m(\mathbb T^3)=H^m$,
  $L^2(\mathbb T^3)=L^2$ and $L^\infty(\setT^3)=L^\infty$.
\end{defn}

\subsection{Local existence}
\label{sec:local-existence}
In this subsection, we address the local existence theorem for the
initial value problem \eqref{M-eq:wave:1}--\eqref{M-eq:wave:2}. 
The theory of second-order quasilinear wave equations of the
type 
\begin{equation} 
  \label{eq:section3-fourier2} 
  \begin{cases}
    & g^{\alpha\beta}(u,u')\partial_\alpha\partial_\beta u=\mathbf{F}(u,u')\\
 &   u(0,\cdot) =u_0, \ \partial_t u(0,\cdot)= u_1
  \end{cases} 
\end{equation}
where the metric $g^{\alpha\beta}(u,u')$ is of Lorentzian signature and
$u'=\partial_\alpha u$, is well documented. 
Notable treatments include \cite[Theorem 6.4.11]{Hormander_1997},
\cite[Theorem 4.1]{Sogge_95}, and
\cite[Theorem~5.1]{Shatah-Struwe-98}. 
\makeatletter \def\plain@equationautorefname{Equation} 
Though \autoref{M-eq:wave:1} is semilinear, these references are nonetheless
pertinent, for the principal challenge—the estimation of the nonlinear
source term $\mathbf{F}(u,u')$—is shared between quasi- and semilinear cases
alike.

\begin{thm}[Local existence of solutions to a nonlinear wave equation]
  \label{thm:Sogge}
  Suppose that $g^{\alpha\beta}$ and $\mathbf{F}$ are $C^\infty$ functions with
  bounded derivatives, that the metric $g^{\alpha\beta}$ is a small
  perturbation of the Minkowski metric,
  and that $\mathbf{F}(0,0)=0$. 
  Let $m >\frac{(n+2)}{2}$. 
  If the initial data $(u_0, u_1)$ belong to
  $H^{m+1}(\setR^n)\times H^m(\setR^n)$, then there exists a time
  $T>0$, depending on the norm of the initial data, for which the
  Cauchy problem \eqref{M-eq:section3-fourier2} admits a unique
  solution satisfying
  \begin{equation}
    \label{eq:section2-preliminaries:3}
    u \in L^\infty([0,T];H^{m+1}(\setR^n))\cap C^{0,1}([0,T];H^m(\setR^n)),
  \end{equation}
where $C^{0,1} $ denotes the spaces of Lipschitz functions.
\end{thm}

We consider the Cauchy problem \eqref{M-eq:wave:1}--\eqref{M-eq:wave:2}, on
the three-dimensional torus $\setT^3$ -with a flat metric,  of the form
\begin{equation}
  \label{eq:source}
  \mathbf{F}(u,u')=- 2\Omega \partial_t u + \mathord{e^{-\varkappa t}} a(t,x)(1+u)^{\mu}
\end{equation}
The local existence result for the Cauchy problem
\eqref{eq:section3-fourier2} extends naturally to the  torus
$\setT^3$, where the  periodicity property of the functions
ensure that integration by parts remains straightforward. 

Though in \autoref{eq:source} one observes that $\mathbf{F}(0,0)\neq 0$, this
is not fatal: the inhomogeneity
$\mathbf{F}(0,0)=\mathord{e^{-\varkappa t}} a(t,x)$ lies in $H^m(\setT^3)$, and
the structure of the theorem accommodates such data, provided the
necessary smoothness are retained. 

The obstacle arises with the nonlinear term $(1+u)^\mu$, whose
derivatives are no longer uniformly bounded as $u$ approaches $-1$. 
This nonlinearity imposes a natural constraint on the size of
admissible initial data. 
Nonetheless, if $u_0$ is sufficiently small in norm, then 
\fullref{thm:Sogge} applies to the Cauchy problem
\eqref{M-eq:wave:1}--\eqref{M-eq:wave:2}, ensuring local existence.

\begin{thm}[Local existence of solutions to the Nordström equation]
  \label{thr:sec2-LocalExistence}
  Let $m\geq 3$ and suppose $a(t,\cdot)\in H^m(\setT^3)$. 
  If the initial data $(u_0,u_1)$ belong to
  $H^{m+1}(\setT^3)\times H^m(\setT^3)$ and the norm
  $\|u_0\|_{H^m(\setT^3)}$ is sufficiently small, then there exists a
  positive time $T$ and a unique solution $u$ to the Cauchy problem
  \eqref{M-eq:wave:1}--\eqref{M-eq:wave:2} satisfying
  \begin{equation}
    \label{eq:existence}
    u \in L^\infty([0,T];H^{m+1}(\setT^3))\cap C^{0,1}([0,T];H^m(\setT^3)).
  \end{equation}
\end{thm}

\begin{proof}[~\autoref{thr:sec2-LocalExistence}] 

If $\|u_0\|_{H^m} $ is small, then by continuity $\|u(t)\|_{H^{m}}$ 
remains small in a certain time interval $[0,t_1)$, then $\mathbf{F}(u(t),u'(t))$ is a smooth 
function and we can apply Theorem \fullref{thm:Sogge}  to the Cauchy
  problem \eqref{M-eq:wave:1}--\eqref{M-eq:wave:2}, thereby establishing
  the result. 

\end{proof}

\subsection{Calculus in the Sobolev spaces on the torus $\setT^3$}
\label{sec:calculus-sobol}
 
We summarise here several established results concerning calculus in
Sobolev spaces. 
The reader is referred to \cite[Chap.~6.4]{Hormander_1997} for detailed
proofs. 
Although the original statements are framed in $\setR^n$, the arguments
extend naturally to the torus $\setT^n$ (see also \cite{taylor97}). 
While the literature typically presents these results in the broader
context of Sobolev spaces $W^{m,p}$, we restrict attention to the
Hilbertian case $H^m(\setT^n)$, which suffices for our purposes.
\begin{prop}[Sobolev inequality
  \texorpdfstring{\protect\cite[Cor.~6.4.9]{Hormander_1997}}{}]
  \label{prop:Sobolev-inequalyi}
  If $m$ is a positive integer, $\frac{n}{ 2}<m$ and $u\in H^m(\setT^n)$, then $u$ is a
  continuous function and
  \begin{equation}
    \label{eq:Sobolev:1}
    \|u\|_{L^\infty} \leq C_m \|u\|_{H^m}.
  \end{equation}
\end{prop}

\begin{prop}[Multiplications of derivatives in $L^2(\setT^n)$
  \texorpdfstring{\protect\cite[Cor.~6.4.4]{Hormander_1997}}{}]
  \label{prop:multiplication}
  If $u,v\in L^\infty(\setT^n)$,
  $\partial_{\mathord{\vec{\alpha}}} u, \partial_{\mathord{\vec{\alpha}}} v \in 
L^2(\setT^n)$, $|{\mathord{\vec{\alpha}}}|=k$ and $m$ is a 
nonnegative integer, then
  $\partial_{\mathord{\vec{\alpha}}}(uv)\in L^2(\setT^n)$ when 
$|{\mathord{\vec{\alpha}}}|=k$ and
  \begin{equation}
    \|\partial_{\mathord{\vec{\alpha}}}(uv)\|_{L^2}\leq 
    C_k\left(\|v\|_{L^\infty}\sum_{|{\mathord{\vec{\alpha}}}|=k}\|\partial_{\mathord{\vec{\alpha}}} u\|_{L^2} + 
      \|u\|_{L^\infty}\sum_{|{\mathord{\vec{\alpha}}}|=k}\|\partial_{\mathord{\vec{\alpha}}} v\|_{L^2}\right).
  \end{equation}
\end{prop}

Combining  \fullref{prop:Sobolev-inequalyi} and
\fullref{prop:multiplication}, we obtain
\begin{cor}[Algebra of $H^m(\setT^n)$]
  \label{cor:algebra}
  If $u,v\in H^m$ and $m > n/2$, then $uv\in H^m(\setT^n)$ and
  \begin{equation}
    \|uv\|_{H^m}\leq C_m 
    \left(\|v\|_{L^\infty}\|u\|_{H^m}+\|u\|_{L^\infty}\|v\|_{H^m}\right)\leq C_m \|u\|_{H^m}\|v\|_{H^m}.
  \end{equation}
\end{cor}

\begin{prop}[Moser type estimate
  \texorpdfstring{\protect\cite[Cor.~6.4.5]{Hormander_1997}}{}]
  \label{prop:Moser}
  Let $u\in L^\infty(\setT^n,\setR)$, $F\in C^k(\setR) $ and assume
  $\partial_{\mathord{\vec{\alpha}}} u \in L^2(\setT^n)$. Then
  $\partial_{\mathord{\vec{\alpha}}} \left(F(u)\right)\in L^2(\setT^n)$ when
  $|{\mathord{\vec{\alpha}}}|=k$, and
  \begin{equation}
    \label{eq:Moser:1}
    \sup_{|{\mathord{\vec{\alpha}}}|=k}\|\partial_{\mathord{\vec{\alpha}}}\left( 
      F(u)\right)\|_{L^2}\leq C_k\sup_{1\leq l\leq 
      k}|F^{(l)}(u)|\ 
    \left(\|u\|_{L^\infty}\right)^{l-1}\ \sup_{|{\mathord{\vec{\alpha}}}|=k}\|\partial_{\mathord{\vec{\alpha}}} u\|_{L^2}, 
  \end{equation}
  when $k>0$ and for $k=0$
  \begin{equation}
    \label{eq:Moser:2}
    \| F(u)-F(0)\|_{L^2}\leq M \|u\|_{L^2},
  \end{equation}
  where $M$ is the Lipschitz constant of $F$.
\end{prop}

 \begin{prop}[Fractional power estimate]
   \label{prop:fractional}
   Let $\mu\in\setR$. 
   If $u\in H^m(\setT^n)$, $m >n/2$ and $\|u\|_{L^\infty}\leq \delta'<1$, then there is a
   positive constant $C_{m,\mu,\delta'}$ such that
   \begin{equation}
     \label{eq:Moser:3}
     \|(1+u)^\mu\|_{H^m}\leq  C_{m,\mu,\delta'} \|u\|_{H^m}+(2\pi)^{\frac{3}{2}}.
   \end{equation}
 \end{prop}

 \begin{rem}
   We can replace the condition $\|u\|_{L^\infty}\leq \delta'$ by
   $\|u\|_{H^m}\leq \delta$ for some positive $\delta$.
   \makeatletter
   \def\plain@equationautorefname{Inequality}
   This follows from from Sobolev's  \autoref{eq:Sobolev:1}.
 \end{rem}
 
 \begin{proof}[~\autoref{prop:fractional}]
   We invoke Moser’s result (\autoref{prop:Moser}) with the function
   $ F(x) = (1 + x)^\mu $, under the assumption that $ |x|\leq \delta' $ holds. 
   Then for any positive integer $ l$ we obtain
   \begin{equation}
     F^{(l)}(x)=\mu(\mu-1)\cdot (\mu-(l-1))(1+x)^{\mu-l}.
   \end{equation}
   So for $|x|\leq \delta'$, we obtain the following inequality
   \begin{equation}
     |F^{(l)}(x)|\leq |\mu(\mu-1)\cdot (\mu-(l-1))|
     \left\{\begin{array}{ll}
       \left(1+\delta'\right)^{\mu-l}, & \mu-l\geq 0 \\ 
       (1-\delta')^{-(l-\mu)}, & \mu-l<0
     \end{array}
   \right\}
   \upperrel{\rm def}{=}c_l.
 \end{equation}
 Note that $\left(\|u\|_{L^\infty}\right)^{l-1}\leq \left(1+\delta'\right)^{l-1}$ holds
 Let $ M_{k,\delta'}=\max_{1\leq l\leq k}c_l \left(1+\delta'\right)^{l-1}$, then
   \makeatletter
   \def\plain@equationautorefname{Inequality}
 \autoref{eq:Moser:1} implies that for $|{\mathord{\vec{\alpha}}}|=k$,
 \begin{equation}
   \label{eq:Moser:4}
   \|\partial_{\mathord{\vec{\alpha}}}(1+u)^\mu\|_{L^2}\leq C_kM_{k,\delta'} \sup_{|{\mathord{\vec{\alpha}}}|=k} 
   \|\partial_{\mathord{\vec{\alpha}}} u\|_{L^2}
 \end{equation}
 holds, if $k>0$. 
 For $k=0$, the Lipschitz constant $M$ is given by 
 \begin{equation*}
   M=M_{\mu,\delta'}=\sup_{|x|\leq \delta'} |\mu| (1+x)^{\mu-1},
 \end{equation*}
 hence, \autoref{eq:Moser:2}
 implies that
 \begin{equation}
   \label{eq:Moser:5}
   \|(1+u)^\mu\|_{L^2}\leq M_{\mu,\delta'}  \|u\|_{L^2}+\| 1\|_{L^2}=M_{\mu,\delta'}  
   \|u\|_{L^2}+(2\pi)^{\frac{3}{2}}
 \end{equation}
 holds. Let
 \begin{equation*}
   C_{m,\mu}=\max_{1\leq k\leq m}\{C_k,M_{k,\delta'} M_{\mu,\delta'}\},
 \end{equation*}
 then \autoref{eq:Moser:4} and \autoref{eq:Moser:5} imply
 \autoref{eq:Moser:3}.

\end{proof}

\subsection{Gronwall inequality}
We shall use the following version of Gronwall's inequality (see 
e.~g.~\cite{Bahouri_2011}).
\begin{lem}[Gronwall's inequality]
  \label{lem:section2-preliminaries:1}
  Let $ f $ and $ g $ be nonnegative functions on $[t_0, T]$, with $ f $
  continuous and $ g $ differentiable.
  Let $ A $ be a continuous function on $[t_0, T]$.
  Suppose that for $t\in [t_0,T]$,
  \begin{equation}
    \label{eq:section2-preliminaries:1}
    \frac{1}{2}\frac{d}{dt}g^2(t)\leq A(t)g^2(t)+f(t)g(t),
  \end{equation}
  holds, then for $t\in [t_0,T]$ we have the following inequality 
  \begin{equation}
    \label{eq:section2-preliminaries:2}
    g(t)\leq e^{\int_{t_0}^t A(\tau)d\tau} g(t_0)+\int_{t_0}^t e^{\int_{\tau}^tA(s)ds} f(\tau)d\tau.
  \end{equation}
\end{lem}
\makeatletter 
  \def\plain@equationautorefname{Equation}

\renewenvironment{proof}[1][\proofname]{{\bfseries Proof of #1.}}{\hfill$\blacksquare$}
\makeatletter 
  \def\plain@equationautorefname{Equation}

\section{Energy estimates and decaying estimates}
\label{sec:energy-estimates}

In this section, we introduce our primary tool: the energy estimates, 
which enable us to state and prove
\fullref{M-lem:section4-main-result:1} and
\fullref{M-lem:section4-main-result:2}, which then in turn allow us to
prove the global existence result.

We will derive the energy estimates for \autoref{M-eq:wave:1}, 
which we write as
\begin{equation}
  \label{eq:linear-wave}
  \partial_t^2u-\Delta u =-2\Omega \partial_t u +F \qquad F(t,x,u):= \mathord{e^{-\varkappa t}} a(t,x)(1+u)^{\mu}  
\end{equation}
\begin{rem}
  \label{rem:section3-energy-estimates:4}
  Please note that the Nordström equation is a semi-linear wave
  equation, and that the nonlinear nature of the lower-order
  terms—denoted by $F$—does not play a significant role in the energy
  estimates. 
  In that sense, \autoref{eq:linear-wave} may be regarded as a linear
  wave equation.

  This observation simplifies the derivation of higher-order energy
  estimates, as there are no commutator terms to contend with, unlike
  in the case of quasilinear wave equations. 
  For this reason, we present the $L^2$ estimate in detail, from
  which the higher-order estimates follow directly.
 
\end{rem}

So we start with the following definition.

\begin{defn}[Higher order energy]
  \label{def:decaying-energy-estimates-nordstroem:2}
  We define the \emph{positive definite} energies $\mathcal{E}$ and
  $ E_{m}$ as follows:
  \begin{equation}
    \label{eq:L-2-energy}
    \mathcal{E}^2[u]\upperrel{\rm def}{=} \frac{1}{2} \int_{\mathbb{T}^3} \left[      (\partial_t u)^2 + \Omega  u  \partial_{t}u +\frac{1}{2} \Omega^{2}u^{2} \right] d^{3}x +
    \frac{1}{2}\int_{\mathbb{T}^3} |\nabla u|^2 d^{3}x,
  \end{equation}
  where $|\nabla u|^2$ is defined by \autoref{M-eq:section1-intro:14}.
  Let $\mathord{\vec{\alpha}}$ be a spatial multi-index, then
  $ E_{m}^2 $ is defined as
  \begin{equation}
    \label{eq:decaying-energy-estimates-nordstroem:3A}
    E_{m}^2  \overset{\mbox{\tiny{def}}}{=} \sum_{0 \leq |\vec{\alpha}| \leq m} \mathcal{E}^2[ \partial_{\mathord{\vec{\alpha}}} u],
  \end{equation}

  $E_{m}^{2}$ can also be written as:
  \begin{equation}
    \label{eq:decaying-energy-estimates-nordstroem:3}
    \begin{aligned}
      E_{m}^2  &=
               \sum_{0 \leq |\vec{\alpha}| \leq m}
               \frac{1}{2} \int_{\mathbb{T}^3} 
               (\partial_t \partial_{\vec{\alpha}}u)^2
               + \Omega \left(  \partial_{\vec{\alpha}} u \right)\left(\partial_{\vec{\alpha}} \partial_{t}u \right)
               +\frac{1}{2} \left(\Omega\right)^{2}\left(\partial_{\vec{\alpha}} u\right)^{2} \,d^3x\\
             &  +
               \sum_{0 \leq |\vec{\alpha}| \leq m} \frac{1}{2} \int_{\mathbb{T}^3} 
               |\nabla (\partial_{\vec{\alpha}}u)|^2 d^3x
    \end{aligned}
  \end{equation}
  
\end{defn}

\begin{rem}[About the Choice of the Energy Functional]
  \label{rem:section3-energy-estimates:3}
  This particular choice of coefficients in \autoref{eq:L-2-energy}
  satisfies two important properties of the energy $\mathcal{E}^2$.
\begin{enumerate}
  \item 
  First, the quadratic form
  \begin{equation}
    \label{eq:energy:2}
    \left( \partial_t u, u\right)\upperrel{\rm def}{=} \int_{\mathbb{T}^3} 
    \left[ (\partial_t u)^2
      + \Omega u \partial_t u 
      + \frac{1}{2} \Omega^2 u^2 \right] d^3x
  \end{equation} 
  is positive definite.
  
  \item 
  Second, we address the indefinite quadratic form
  $\left( \partial_t u, u \right) = \int_{\mathbb{T}^3} \partial_t u \, u \, d^3x$,
  which appears in both $\mathcal{E}^2$ and its time derivative.
  Since this form is indefinite, it precludes the direct application of standard inequalities.
  Nonetheless, the specific choice of coefficients permits us to absorb this term
  into $-\Omega \mathcal{E}^2$, thereby retaining control over the energy.
\end{enumerate}
\end{rem}

We begin by stating and proving the essential $L^{2}$ estimate. 
The semi-linear nature of \autoref{eq:linear-wave} then permits a straightforward generalisation
to higher-order energy estimates.

\begin{prop}[The $L^{2}$  energy estimate for the Nordström equation]
  \label{prop:energy:2}
  Let $ u $ be a solution to equation \eqref{eq:linear-wave}, then the following inequality holds:
  \begin{equation}
    \label{eq:section3-energy-estimates:13}
    \frac{d}{dt} \left(\mathcal{E}^2[u]\right) (t)\leq -\Omega \mathcal{E}^2[u] (t)+\left(\frac{\Omega^2}{\sqrt{2}}\|u(t)\|_{L^2} 
    +\sqrt{2} \|F\|_{L^2}\right)\mathcal{E}[u](t).
  \end{equation}
\end{prop}

The proof of that proposition is done with the help of the following lemma.
\begin{lem}[An intermediate $L^{2}$ estimate]
  \label{prop:energy:1}
  Let $ u $ be a solution to \autoref{eq:linear-wave}, then the
  following inequality
  \begin{equation}
    \label{eq:energy:3}
    \frac{d}{dt} \mathcal{E}^2[u]\leq -\Omega \mathcal{E}^2[u]+\frac{\Omega^3}{4}\|u\|^2_{L^2}+\int_{\setT^3}\left(\partial_t u+\frac{\Omega}{2} u\right)F d^{3}x.
  \end{equation}
  is true.
\end{lem}
\begin{proof}[~\autoref{prop:energy:1}]
    We start with differentiating the energy with respect to $t$
    \begin{equation}
      \label{eq:section3-energy-estimates:10}
      \begin{aligned}
        \frac{d}{dt}\mathcal{E}^2[u] &  =
                            \frac{1}{2}\frac{d}{dt} \int_{\mathbb{T}^3} \left[(\partial_t u)^2  +  |\nabla u|^2+\Omega  u  \partial_{t}u +\frac{1}{2} \Omega^{2}u^{2} \right] d^{3}x \\
                          &=
                            \int_{\mathbb{T}^3}\left[ ( \partial_t u) (\partial_t^2 u)+\left(\nabla u\cdot \partial_t(\nabla u)\right)
                            +\frac{\Omega}{2} \left(u \partial_t^2 u\right) + \frac{\Omega}{2} (\partial_t u)^2 + \frac{\Omega^2}{2} u \partial_t u \right]d^{3}x
      \end{aligned}
    \end{equation}
    Integration by parts leads to
    \begin{align}
      \label{eq:proof:2}
      \frac{d}{dt}\mathcal{E}^2[u] &  =
                          \int_{\mathbb{T}^3}\bigg[\partial_t u \partial_t^2 u + \underbrace{\nabla u\cdot \partial_t(\nabla u)}_{=-( \partial_t u)\Delta u}    
                          +\frac{\Omega}{2} \left(u \partial_t^2 u\right) + \frac{\Omega}{2} (\partial_t u)^2 +\frac{\Omega^2}{2} u \partial_t u \bigg]\,d^{3}x \\
      \label{eq:proof:1}
                        &=    \int_{\mathbb{T}^3}\left[\partial_t u 
                          \left( \partial_t^2 u -\Delta u \right)
                          +\frac{\Omega}{2} \left(u \partial_t^2 u\right)
                          + \frac{\Omega}{2} (\partial_t u)^2 +\frac{\Omega^2}{2} u \partial_t u \right]\,d^{3}x \\
    \end{align}
    Inserting \autoref{eq:linear-wave} into \autoref{eq:proof:1} leads
    to
    \begin{equation}
      \label{eq:proof:3}
      \frac{d}{dt}\mathcal{E}^2[u]   =
      \int_{\mathbb{T}^3}\left[\partial_t u \left( -2\Omega \partial_t u +F \right)
        +\frac{\Omega}{2}u \left(\Delta u -2\Omega \partial_t u +F  \right)
        + \frac{\Omega}{2} (\partial_t u)^2 +\frac{\Omega^2}{2} u \partial_t u \right]\,d^{3}x 
    \end{equation}
    Again integration by parts of the term
    $\int_{\mathbb{T}^3}\frac{\Omega}{2}u \Delta u d^{3}x$ in \autoref{eq:proof:3} leads to
    \begin{equation}
      \label{eq:proof:4}
      \frac{d}{dt}\mathcal{E}^2[u]   =
      \int_{\mathbb{T}^3}\left[\partial_t u \left( -2\Omega \partial_t u +F \right)
        -\frac{\Omega}{2} \left|\nabla u  \right|^{2}  -\Omega^{2}u \partial_t u + \frac{\Omega}{2}u F 
        + \frac{\Omega}{2} (\partial_t u)^2 +\frac{\Omega^2}{2} u \partial_t u \right]\,d^{3}x 
    \end{equation}
    Simplifying terms in \autoref{eq:proof:4} results in
    \begin{align}
      \label{eq:proof:5}
      \frac{d}{dt}\mathcal{E}^2[u]
      &=
        \int_{\mathbb{T}^3}\left[ \left( -2\Omega+\frac{\Omega}{2}\right ) (\partial_t u)^{2} + \partial_t uF
        -\frac{\Omega}{2} \left|\nabla u  \right|^{2}+  \left( - \Omega^{2}u + \frac{\Omega^{2}}{2} \right) \partial_t u
        + \frac{\Omega}{2}u F \right]\,d^{3}x                    \\
      &=
        \int_{\mathbb{T}^3}\left[ -\frac{3}{2}\Omega (\partial_t u)^{2}  + \partial_{t}uF
        -\frac{\Omega}{2} \left|\nabla u  \right|^{2}  - \frac{1}{2}\Omega^{2}u \partial_t u + \frac{\Omega}{2}u F
        + \frac{\Omega}{2} (\partial_t u)^2  \right]\,d^{3}x              \\
      & =
        -\frac{\Omega}{2}
        \int_{\mathbb{T}^3} \left[3(\partial_t u)^2+|\nabla u|^2 +\Omega u (\partial_t u)\right]
        d^{3}x + \int_{\mathbb{T}^3}\left(\partial_t u+\frac{\Omega}{2}u\right)F d^{3}x \\
      &  \leq
        -\frac{\Omega}{2}
        \int_{\mathbb{T}^3} \left[(\partial_t u)^2+|\nabla u|^2 +\Omega u (\partial_t u)\right]
        d^{3}x + \int_{\mathbb{T}^3}\left(\partial_t u+\frac{\Omega}{2}u\right)F d^{3}x \\
      & = -\Omega
        \mathcal{E}^2[u]+\frac{\Omega^3}{4}\|u\|_{L^2}^2
        + \int_{\mathbb{T}^3}\left(\partial_t u+\frac{\Omega}{2}u\right)F d^{3}x.
    \end{align}
    which finishes the proof.
  \end{proof}

\makeatletter 
  \def\plain@equationautorefname{Inequality}

We turn now to the proof \fullref{prop:energy:2}.\\
\begin{proof}[~\autoref{prop:energy:2}]
  We need to estimate the terms $\| u\|_{L^2} $ and
  $ \|\partial_t u+ \frac{\Omega}{2}u\|_{L^2}$ by the energy $\mathcal{E}[u]$. 
  \begin{enumerate}
    \item 
    So we first claim that
    \begin{equation}
      \label{eq:L-2-energy:2}
      \| u\|_{L^2} \leq \frac{\sqrt{8}}{\Omega} \mathcal{E}[u].
    \end{equation}
    holds. Taking into account the definition of the energy (namely~\autoref{eq:L-2-energy}), it suffices to show that
    \begin{equation}
      \| u\|_2^2\leq \frac{8}{\Omega^2}  \frac{1}{2} \int_{\mathbb{T}^3} \left[      (\partial_t u)^2 + \Omega  u  \partial_{t}u +\frac{1}{2} \Omega^{2}u^{2} \right] d^{3}x,
    \end{equation}
    is satisfied, which can be written as 
    \begin{equation}
      0\leq \frac{8}{\Omega^2}  \frac{1}{2}\int_{\mathbb{T}^3}\left[ ( \partial_t u)^2 +\Omega \left(u \partial_t  u\right)u+{\Omega^2}  u^2\left(\frac{1}{2}-\frac{1}{4}\right)\right] d^{3}x.
    \end{equation}
    Observing that
    \begin{equation}
      \int_{\mathbb{T}^3}\left[ ( \partial_t u)^2 +\Omega \left(u \partial_t  u\right)u+\frac{{\Omega^2}  u^2}{4}\right] d^{3}x =\int_{\mathbb{T}^3}\left[  \partial_t u +\frac{\Omega}{2}u\right]^2 d^{3}x,
    \end{equation}
    is true, we see that the following inequality 
    \begin{equation*}
      \| u\|_{L^2}^{2} \leq \frac{{8}}{\Omega^{2}} \mathcal{E}^{2}[u].
    \end{equation*}
    holds, which after taking the square roots, is the desired  \autoref{eq:L-2-energy:2}.
    \item 
    It is simpler to estimate the second term $\left \|\partial_t u+ \frac{\Omega}{2}u\right\|_{L^2}$ by the energy $\mathcal{E}[u]$.
    We  considering the following inequalities
    \begin{equation}
      \label{eq:energy-L-2:4}
      \begin{aligned}
        \left\|\partial_t u+ \frac{\Omega}{2}u\right\|_{L^2}^2
        & = \int_{\setT^3}\left[ (\partial_t u)^2 +\Omega \left(u \partial_t  u\right)u+\frac{\Omega^2}{4}  u^2\right] d^{3}x \\
        &\leq 
          \int_{\setT^3}\left[ ( \partial_t u)^2 +\Omega \left(u \partial_t  u\right)u+\frac{\Omega^2}{2}  u^2\right]d^{3}x \leq2  \mathcal{E}^2[u].
      \end{aligned}
    \end{equation}

    We observe that by the Cauchy Schwarz inequality and \autoref{eq:energy-L-2:4}, we obtain the following inequality
    \begin{equation}
      \label{eq:section3-energy-estimates:8}
      \int_{\mathbb{T}^3}\left(\partial_t u+\frac{\Omega}{2}u\right)F d^{3}x
      \leq \left\|\partial_t u+\frac{\Omega}{2}u\right\|_{L^2} \|F\|_{L^2}\leq \sqrt{2}\mathcal{E}[u] \|F\|_{L^2}.
    \end{equation}
    Combining \autoref{eq:section3-energy-estimates:8},
    \autoref{eq:energy:3}, the Cauchy-Schwarz inequality and
    \autoref{eq:L-2-energy:2} results in the following inequality  
    \begin{equation}
      \label{eq:section3-energy-estimates:11}
      \begin{aligned}
        \frac{d}{dt} \mathcal{E}^2[u] &\leq -\Omega \mathcal{E}^2[u]+\frac{\Omega^3}{4}\|u\|^2_{L^2}+\int_{\setT^3}\left(\partial_t u+\frac{\Omega}{2} u\right)F d^{3}x. \\
                           &\leq-\Omega \mathcal{E}^2[u] + \frac{\Omega^3}{4}\|u\|_{L^2} \frac{\sqrt{8}}{\Omega} \mathcal{E}[u] + \sqrt{2}\mathcal{E}[u] \|F\|_{L^2} \\
                           &=-\Omega \mathcal{E}^2[u] + \frac{\Omega^2}{\sqrt{2}}\|u\|_{L^2}  \mathcal{E}[u] + \sqrt{2}\mathcal{E}[u] \|F\|_{L^2}                 \\
      \end{aligned}
    \end{equation}
    which finishes the proof.
  \end{enumerate}
\end{proof}

\makeatletter 
  \def\plain@equationautorefname{Equation}

  \begin{prop}[Higher energy estimates]
    \label{prop:decaying-energy-estimates-nordstroem:2}
    Let $u$ be a solution to \autoref{eq:linear-wave}, and $E_{m}^{2}$
    its corresponding energy defined by
    \autoref{eq:decaying-energy-estimates-nordstroem:3}, then the
    following inequality
    \begin{equation}
      \label{eq:decaying-energy-estimates-nordstroem:5}
      \frac{d E_m^{2}}{dt}(t) \leq -\Omega 
      E_{m}^{2}(t)+\left(\frac{\Omega^{2}}{\sqrt{2}}\|u\|_{H^{m}(\setT^{3})} + 
        \sqrt{2}\|F\|_{H^m(\setT^{3})}\right)E_{m}(t).
    \end{equation}
    \noeqref{eq:decaying-energy-estimates-nordstroem:5} holds.
  \end{prop}
  
\makeatletter 
  \def\plain@equationautorefname{Inequality}

\begin{proof}[~\autoref{prop:decaying-energy-estimates-nordstroem:2}]
  Since $\partial_{\mathord{\vec{\alpha}}}u$ satisfies the following equation
  \begin{equation}
    \label{eq:linear-wave:2}
    \partial_t^2(\partial_{\mathord{\vec{\alpha}}}u)-\Delta(\partial_{\mathord{\vec{\alpha}}}u)  =-2\Omega 
    \partial_t(\partial_{\mathord{\vec{\alpha}}}u)  +\partial_{\mathord{\vec{\alpha}}}F,
  \end{equation}
  we can apply  \fullref{prop:energy:1} and
  \fullref{prop:energy:2} to each term of the energy end obtain
  \autoref{eq:decaying-energy-estimates-nordstroem:5}

\end{proof}

\begin{rem}[The Strategy for Proving Global Existence]
  Our approach to establish global existence hinges on maintaining
  a small ratio $\frac{\|u\|_{H^m} }{ \|\partial_t u\|_{H^m}}$ for all $t \geq 0$. 
  This condition ensures that the term $\|u\|_{H^m}$ does not obstruct the
  decay of the energy $E_m(t)$.

  The main result of this work is subject to the restriction
  $\Omega <1$. 
  Recall that $2\Omega$ is the coefficient of the dissipative term $\partial_t u$. 
  At first glance, this constraint may seem counterintuitive, since
  stronger dissipation typically enhances decay.

  The rationale for imposing this restriction lies in the structure of
  the energy estimate
  \eqref{eq:decaying-energy-estimates-nordstroem:5}, which contains
  the term $\|u\|_{H^m}$. 
  If the dissipative term $2\Omega \partial_t u$ becomes too large, it risks
  undermining the smallness of the ratio
  $\frac{\| u\|_{H^m} }{ \|\partial_t u\|_{H^m}}$, thereby invalidating a key
  assumption in the analysis.
\end{rem}

We now state the following corollary, which follows easily from
\autoref{prop:decaying-energy-estimates-nordstroem:2} and
\fullref{M-lem:section2-preliminaries:1}

\begin{cor}[Final energy estimates for the gravitational field]
  \label{cor:section3-energy-estimates:1}
  Under the
  hypotheses of \fullref{prop:decaying-energy-estimates-nordstroem:2}, the following inequality 
  \begin{equation}
    \label{eq:section3-energy-estimates:9}
    \begin{aligned}
      E_{m}(t)
      & \leq e^{-\Omega (t-t_0)}E_{m}(t_0)
        +\frac{\Omega^2}{\sqrt{2}}\int_{t_0}^te^{-\Omega (t-s)} \left\Vert u(s) \right\Vert_{H^{m}}ds 
        + \int_{t_0}^te^{-\Omega  (t-s)} \sqrt{2} \left\Vert F(s) \right\Vert_{H^{m}}ds
    \end{aligned}.
  \end{equation}
holds
\end{cor}

\renewenvironment{proof}[1][\proofname]{{\bfseries Proof of #1.}}{\hfill$\blacksquare$}
\makeatletter
  \def\plain@equationautorefname{Equation}
\section{Main results}
\label{sec:main-results}
\begin{thm}[Global existence]
  \label{thr:section4-main-result:1}
  Let
  \begin{equation}
    \label{eq:section4-main-result:1}
    0<  \Omega<1 \quad 3\leq m, \quad \mu\in \setR, \quad a(\cdot,t)\in H^m(\setT^{3}), \quad \forall  t\geq 0.
  \end{equation}
  Let the initial data of \autoref{M-eq:wave:1}
  $u_0\in H^{m+1}(\setT^{3})$ and $ u_1\in H^m(\setT^{3})$. We assume that these initial
  data have zero mean and that they satisfy the inequality
  \begin{equation}
    \label{eq:section4-main-result:8}
    \frac{1}{4}\|u_0\|_{H^m} \leq E_m(0)
  \end{equation}
  (we can always choose such initial data). We also assume the
  following inequality
  \begin{equation}
    \label{eq:section4-main-result:9}
    E_m(0)\leq \delta{\Omega} 
  \end{equation}
  where $\delta$ is a positive constant depending on the Sobolev's
  embedding theorem. Moreover, we assume, that there exists a
  $\varepsilon_{0}>0$, such that for every $\varepsilon$ with $0<\varepsilon\leq\varepsilon_{0}$,
  \begin{equation}
    \label{eq:small:1}
    \sup_{[0,\infty)} \left\Vert a(t,\cdot) \right\Vert_{H^{m}} \leq \varepsilon
  \end{equation}
  holds. Then the classical solution $u$ of the Cauchy problem
  \eqref{M-eq:wave:1}-\eqref{M-eq:wave:2} exists on
  $[0,\infty) \times \mathbb{T}^3$ and the inequality
  \begin{equation}
    \label{eq:section4-main-result:3}
    \| u(t)\|_{H^{m}}\leq \sqrt{2}\delta \Omega
  \end{equation}
  holds on $[0,\infty) \times \mathbb{T}^3$.
\end{thm}
The proof of this theorem relies on  a standard bootstrap argument.
We now impose the following.
\begin{assum}[Bootstrap assumption]
\label{ass:section4-main-result:1}
Assume the bootstrap condition
\begin{equation}
  \label{eq:section4-main-result:6}
  \frac{1}{2}  \left\Vert u(t) \right\Vert_{H^{m}}^{2}\leq E^{2}_m(0)
\end{equation}
holds.
\end{assum}
\makeatletter 
  \def\plain@equationautorefname{Inequality}
Moreover we make use of the following definition.
\begin{defn}[Maximal time interval]
  \label{def:global-existence-lavi-bootstrap-new:1}
  \begin{equation}
    \label{eq:section4-main-result:2}
    T_{\mathrm{max}} \upperrel{\rm  def}{=} \sup_{T \geq 0}\left\{ 
      \begin{array}{l}     \text{
        The solution $u(t)$ to the
        Cauchy problem \eqref{M-eq:wave:1}–\eqref{M-eq:wave:2} exists}\\
        \text{with the  regularity conditions 
        \eqref{M-eq:existence} for } \ t\in  [0,T),\\
        \text{and   \autoref{eq:section4-main-result:6} 
        is satisfied for} \ t\in  [0,T).
      \end{array}\right.
  \end{equation}
\end{defn}
In order to implement this strategy, a simple to prove proposition and  two essential lemmas are
required, whose proofs are deferred to the appendix
(\autoref{M-sec:proof-fullr-lem:s}). 

\begin{prop}[The bootstrap assumption holds in \texorpdfstring{\protect $[0,T_{1}]$}{}]
\label{prop:section4-main-result:1}
  There exists a $ T_{1}>0$ such that
  \begin{equation}
    \label{eq:section4-main-result:4}
    \| u(t)\|_{H^m}\leq \frac{1}{2}E_m(0) \quad \text{for all}\ \ t\in [0,T_1]
  \end{equation}
  holds. Consequently, $ 0<T_{max}$.
\end{prop}
 
\begin{proof}[~\autoref{prop:section4-main-result:1}]
  By conditions \eqref{M-eq:section4-main-result:9},
  \eqref{M-eq:section4-main-result:8} and \autoref{M-eq:L-2-energy:2},
  we conclude that inequality
  \begin{equation}
    \label{eq:section4-main-result:5}
    \| u_0\|_{H^m}\leq \frac{\sqrt{8}}{\Omega}E_m(0)\leq \sqrt{8} \delta
  \end{equation}
  holds. Hence, by \fullref{M-prop:Sobolev-inequalyi} we conclude
  $ \|u_0\|_{L^\infty}\leq \delta'$ is satisfied for some
  $0<\delta'< 1$. We can now apply \fullref{M-prop:fractional} and conclude
  that $\|(1+u_0)^\mu\|_{H^m}$ is bounded by a constant depending on
  $\delta'$. At that stage we apply \fullref{M-thr:sec2-LocalExistence}
  which tells us that the solution $ u(t)$ exists,
  \begin{equation*}
    u \in L^\infty([0,T_2];H^{m+1}(\setT^3))\cap C^{0,1}([0,T_2];H^m(\setT^3))
  \end{equation*}
  for some $ T_2>0$. Using the initial condition
  \eqref{M-eq:section4-main-result:8}, and the continuity of the norm,
  we conclude that there exists $ 0<T_1\leq T_2 $ such that
  \begin{equation*}
    \|u(t)\|_{H^m}\leq \frac{1}{2}E(0) \qquad \forall  t\in [0,T_1].
  \end{equation*}

\end{proof}

\begin{lem}[An improved energy estimate]
  \label{lem:section4-main-result:1}
  Assume the following:
  \begin{itemize}
    \item Conditions \eqref{eq:section4-main-result:1} of
    \fullref{thr:section4-main-result:1} hold;
    \item the bootstrap assumption is satisfied;
    \item \label{item:global-existence-lavi-bootstrap-new:12}
    \makeatletter \def\plain@equationautorefname{Equation} the
    solution $u(t)$ to the Cauchy problem
    \autoref{M-eq:wave:1}–\autoref{M-eq:wave:2} exists for
    $0\leq t \leq T$, and meets the regularity conditions specified in
    \eqref{M-eq:existence} (\fullref{M-thr:sec2-LocalExistence});
    \makeatletter \def\plain@equationautorefname{Inequality}
    \item There exists $\varepsilon_{1}>0$ such that
    \begin{equation}
      \label{eq:global-existence-lavi-bootstrap-new:70}
      \sup_{[0,\infty)} \left\Vert a(t,\cdot) \right\Vert_{H^{m}} \leq \varepsilon_1
    \end{equation}
  \end{itemize}
  holds, moreover there exists a constant $0<\varepsilon^{\prime}<1$ such that
  \begin{equation}
    \label{eq:global-existence-lavi-bootstrap-new:74}
    E_m(t) \leq \left(1 - \varepsilon'\right) E_m(0) \quad \text{for all } t \in [T_1,T],
  \end{equation}
  is satisfied, provided $\varepsilon_1$ is chosen sufficiently small and the
  existence of $T_1$ is guaranteed by
  \fullref{prop:section4-main-result:1}. The choice of $\varepsilon_1$ depends on
  $\varepsilon'$, which in turn depends on the initial data.
\end{lem}
\begin{lem}[An improved estimate for $\|u(T)\|_{H^{m}}$]
  \label{lem:section4-main-result:2}
  Under the assumptions of 
  \begin{itemize}
    \item \fullref{lem:section4-main-result:1} and
    \item that \autoref{eq:global-existence-lavi-bootstrap-new:74} holds.
  \end{itemize}
  Then then there exists an $\varepsilon_{2}$, satisfying $0< \varepsilon_{2}\leq \varepsilon_{1}$ such that
  \begin{equation}
    \label{eq:strict-inquality:17}
    \frac{1}{2}    \|u(T)\|_{H^{m}}^{2} < E_m^2(0)
  \end{equation}
  holds.
\end{lem}
\subsubsection*{Proof of \fullref{thr:section4-main-result:1}}
\label{sec:proof-fullr-main}
\hfill\\
Now we are in a position to prove our main result.

\begin{proof}[~\protect \autoref{thr:section4-main-result:1}]
  \begin{enumerate}
    \item  There exists
    $ T_{1}>0$ such that
    \begin{equation*}
      \| u(t)\|_{H^m}\leq \frac{1}{2}E_m(0)
    \end{equation*}
    holds for all $t\in [0,T_1]$.
    For the proof we refer to  \textbf{Step 1} in
    the proof of \fullref{lem:section4-main-result:1} in 
    \autoref{M-sec:proof-fullr-lem:s-2}.
    Consequently   $ 0<T_{max}$.
    \item \label{item:global-existence-lavi-bootstrap-new:1}
    Assume  that the maximal time $T_{\mathrm{max}}$
    (as defined in \fullref{def:global-existence-lavi-bootstrap-new:1})
    is finite.
    
    \item Let $\varepsilon_{0}=\min(\varepsilon_{1},\varepsilon_{2})$, where 
    $\varepsilon_{1}$ and $\varepsilon_{2}$ are as in \autoref{lem:section4-main-result:1} 
    and \autoref{lem:section4-main-result:2}.
    Then, for $\varepsilon<\varepsilon_{0}$, 
    \autoref{lem:section4-main-result:1} and \autoref{lem:section4-main-result:2} yield
    \begin{align}
      \label{eq:section4-main-result:7}
      E_m(T_{\mathrm{max}})                          &\leq \left(1 - \varepsilon'\right) 
                                                     E_m(0), 
    \end{align}
    and
    \begin{align}
      \label{eq:section4-main-result:12}
      \frac{1}{2}    \|u(T_{\mathrm{max}})\|_{H^{m}}^{2} &< E_m^2(0).
    \end{align}
    \item \label{item:global-existence-lavi-bootstrap-new:2}
    By continuity of the norm $\|u(t)\|_{H^{m}}^{2}$ in $t$, there exists 
    some $\delta_1 > 0$ such that
    \begin{equation}
      \label{eq:global-existence-lavi-bootstrap-new:60b}
      \frac{1}{2}  \|u(t)\|_{H^{m}}^{2}\leq E^{2}(0) \qquad \text{for all } t\in [T_{\mathrm{max}},T_{\mathrm{max}}+\delta_{1}).
    \end{equation}
    \item
    In light of the bootstrap assumption, it follows that the conditions for the fractional estimate
    \fullref{M-prop:fractional} are satisfied. 
    This permits the application of the
    local existence theorem \fullref{M-thr:sec2-LocalExistence} with the bounded initial data 
    $u(T_{\mathrm{max}})$ and  $\partial_t u(T_{\mathrm{max}})$, which in turn allows us to 
    extend the solution.
    That is, there exists $\delta_2 > 0$ such that the solution remains regular
    and
    \begin{equation}
      \label{eq:global-existence-lavi-bootstrap-new:60}
      u(t) \in C^{0}\left( [0,T_{\mathrm{max}}+\delta_{2});H^{m+1} \right)\cap C^{1}\left( [0,T_{\mathrm{max}}+\delta_{2});H^{m} \right).
    \end{equation}
    \item \label{item:global-existence-lavi-bootstrap-new:4}
    Taken together, these observations imply the existence of a positive time increment 
    $\delta_{3}=\min(\delta_{1},\delta_2)$ such that
    \begin{equation}
      \label{eq:global-existence-lavi-bootstrap-new:63}
      \text{The solution exists classically on } [0,T_{\mathrm{max}}+\delta_{3}) \times \mathbb{T}^3, 
      \quad \text{with }\frac{1}{2} \|u\|_{H^{m}}^{2} \leqslant E^{2}(0).
    \end{equation}
  \end{enumerate}
  This extension beyond $T_{\mathrm{max}}$ contradicts the maximality
  of $T_{\mathrm{max}}$.
  Hence, we are forced to conclude
  \begin{equation}
    \label{eq:global-existence-lavi-bootstrap-new:64}
    T_{\mathrm{max}}=\infty.
  \end{equation}
We also observe that \autoref{eq:global-existence-lavi-bootstrap-new:63} implies that the bootstrap assumption holds 
on $[0,\infty) \times \mathbb{T}^3$ which implies \autoref{eq:section4-main-result:3}.
\end{proof}

\renewenvironment{proof}[1][\proofname]{{\bfseries Proof of #1.}}{\hfill$\blacksquare$}

\section{Proof of \autoref{M-lem:section4-main-result:1} and
  \autoref{M-lem:section4-main-result:2}}
\label{sec:proof-fullr-lem:s}

\subsection{Proof of \fullref{M-lem:section4-main-result:1} }
\label{sec:proof-fullr-lem:s-2}

\hfill
\newline

\begin{proof}[~\protect\fullref{M-lem:section4-main-result:1}]
  \hfill

  \textsc{Proof sketch:} The proof proceeds in three steps: 
  \textbf{1)} Prove the existence of a strictly positive
  $\epsilon^{\prime}$ and an auxiliary function $g(t)$ that is positive for all 
  $  t\geq T_1$. 
  \textbf{2)} Estimate the nonlinear function $F$ by the bootstrap assumptions. 
  \textbf{3)} Use the energy estimates, the nonlinear estimate and the
  results obtained in earlier steps to prove the lemma.

  \renewenvironment{proof}[1][\proofname]{{\bfseries Proof:
    }}{\hfill$\Box$}
  \begin{enumerate}[label=\textbf{Step \arabic*}]

    \item \label{item:section5-appendix:1}
    Let $T_1$ given by \fullref{M-prop:section4-main-result:1} then 
    \begin{itemize}
      \item There exists $0<\varepsilon'<1$ such that
      \begin{equation}
        \label{eq:sec5:20}
        \left[e^{\Omega t}\left(1-\varepsilon'-\Omega\right)-\left(1-\Omega \right)\right]>0, \qquad \text{holds for all} \ t\geq T_1.
      \end{equation}

      \item Set
      \begin{equation}
        \label{eq:sec5:21}
        g(t)=\left[e^{\Omega t}\left(1-\varepsilon'-\Omega\right)-\left(1-\Omega \right)\right]\left[e^{\Omega t}-1\right]^{-1}, 
      \end{equation}
      then
      \begin{equation}
        \label{eq:section5-appendix:3}
        0<g(T_1)\leq g(t),  \qquad \text{for all}\  \ t\geq T_1
      \end{equation}
      holds.
    \end{itemize}

 \begin{rem}
   Note that by the assumption \autoref{M-eq:section4-main-result:1}
   of \fullref{M-thr:section4-main-result:1}, namely $\Omega<1$, we can
   conclude that
   \begin{equation*}
     0<  \left(1-\Omega \right)
   \end{equation*}
   is true, which we will use later.
 \end{rem}

 \begin{proof}
   \autoref{eq:sec5:20} is equivalent to
   \begin{equation}
     \label{eq:sec5:10}
     e^{\Omega t}\varepsilon'<e^{\Omega t}\left(1-\Omega\right) 
     -\left(1-\Omega\right)=\left(e^{\Omega 
         t}-1\right)\left(1-\Omega\right),
   \end{equation}
   or
   \begin{equation}
     \label{eq:sec5:11}
     \varepsilon'<\left(1-e^{-\Omega t}\right)\left(1-\Omega\right)
     \upperrel{\rm def}{=}h(t).
   \end{equation}
   Since $T_1>0$, $h(T_1)>0$, and since
   $\frac{d}{dt}h(t)=\Omega e^{\Omega t} \left(1-\Omega\right)>0$,
   $h(t)\geq h(T_1) $ for all $t\geq T_1$.
   Thus for any $\varepsilon'$ such that satisfies
   \begin{equation}
     \label{eq:sec5:10A}
     \varepsilon'<h(T_1)
   \end{equation}
   \autoref{eq:sec5:20} holds.

   Obviously $g(T_1)>0$.
   To show that $ g(t) $ is an increasing function, we first rewrite
   $ g(t)$:
   \begin{equation}
     \label{eq:sec5:18}
     \begin{aligned}
       g(t) & =\frac{ \left[e^{\Omega 
              t}\left((1-\varepsilon')-\Omega\right) 
              -\left(1-\Omega\right)\right] }{e^{\Omega t}-1}
              = \frac{\left(e^{\Omega t}-1\right)\left(1-\Omega\right)
              -\varepsilon' e^{\Omega t}}{e^{\Omega t}-1} \\
            & = \left(1-\Omega\right)-\frac{\varepsilon' e^{\Omega 
              t}}{e^{\Omega t}-1}
              =\left(1-\Omega\right)-\frac{\varepsilon'}{1-e^{-\Omega t}}.
     \end{aligned}
   \end{equation}
   Hence we conclude the following inequality
   \begin{equation*}
     \frac{d}{dt}g(t) = \frac{ \varepsilon' \Omega e^{-\Omega t}}{(1-e^{-\Omega t})^2}>0
   \end{equation*}
   and the function $g$ is increasing.
 \end{proof}

 \item \label{item:section5-appendix:2}
 \textsc{Claim:} For any $0<T\leq T_{max}$, the following inequality 
 \begin{equation}
   \label{eq:proof:1A}
   \sup_{[0,T]}\|F(u(s))\|_{H^m}  \leq C(\delta^{\prime}) \sup_{[0,\infty)}\|a(s,\cdot)\|_{H^m}.
 \end{equation} 
 is true.
 \begin{proof}
   We start with the function $F(u)$ that is given by
   \begin{equation}
     \label{eq:sec5:12}
     F(u(s))=e^{-\kappa s} a(s, x) (1+u(s))^\mu, \qquad \kappa >0, \mu\in \setR.
   \end{equation}
   By the bootstrap assumption \autoref{M-eq:section4-main-result:6}
   and the condition \eqref{M-eq:section4-main-result:9}, we conclude
   that
   \begin{equation}
     \label{eq:section5-appendix:17}
     \|u(s)\|_{H^m}\leq E_m(0)\leq \Omega \delta, \qquad \text{for all}\ \  s\in [0,t]
   \end{equation}
   holds. 
   Hence we choose $\delta^{\prime}>0 $ so that
   \begin{equation}
     \label{eq:section5-appendix:16}
     \| u(s)\|_{L^\infty}\leq C \|u(s)\|_{H^m}\leq \delta'<1.
   \end{equation}
   is satisfied.
   We now apply \autoref{M-eq:Moser:3} of \fullref{M-prop:fractional}
   together with \autoref{eq:section5-appendix:17} and
   \autoref{eq:section5-appendix:16} to \makeatletter
   \def\plain@equationautorefname{Equation} \autoref{eq:sec5:12} and
   obtain \makeatletter \def\plain@equationautorefname{Inequality} the
   following inequality
   \begin{equation}
     \label{eq:sec5:13}
     \|(1+u(s))^\mu\|_{H^m}\leq C(\delta).
   \end{equation}
   Note that this bound does not depend on $ t $ as long the bootstrap
   assumption \eqref{M-eq:section4-main-result:6} holds and we have
   included, for convenience all constants in $\delta$.
   By \fullref{M-cor:algebra} we finally obtain
   \begin{equation}
     \label{eq:sec5:14}
     \begin{aligned}
       \sup_{[0,t]}\|F(u(s))\|_{H^m} &  \leq \sup_{[0,t]} \|a(s,\cdot)(1+u(s))^\mu \|_{H^m}   
       \\
                                & \leq  C \sup_{[0,t]} \|a(s,\cdot)\|_{H^m} \|(1+u(s))^\mu \|_{H^m} \\
                                & \leq C(\delta) \sup_{[0,t]}\|a(s,\cdot)\|_{H^m} \leq C(\delta) \sup_{[0,\infty)}\|a(s,\cdot)\|_{H^m}, 
     \end{aligned}
   \end{equation}
   where the constant $ C(\delta)$ depends on $\delta$ of condition
   \eqref{M-eq:section4-main-result:9}.
 \end{proof}

 \item \label{item:proof:1}
 \textsc{Claim:} Suppose there exists a strictly positive $\varepsilon_{1}$ such
 that \autoref{M-eq:global-existence-lavi-bootstrap-new:70} is satisfied.
 Then for $\varepsilon^{\prime}$, which is given by inequality \eqref{eq:sec5:10A}  we conclude that 
 \begin{equation}
   \label{eq:section5-appendix:6} 
   E_m(t) \leq \left(1 - \varepsilon'\right) E_m(0) \quad \text{holds for all } t \in 
   [T_1,T].
 \end{equation}

\begin{proof}
  The energy estimate \fullref{M-cor:section3-energy-estimates:1} on
  $[0,t]$ reads
  \begin{equation}
    \label{eq:proof:9}
    \begin{aligned}
      E_m(t) & \leq e^{-\Omega t} E_m(0) + \frac{\Omega^2}{\sqrt{2}}\int_0^t e^{-\Omega(t-s)}\|u(s)\|_{H^m}ds \\
            &\quad+ \sqrt{2}\int_0^t e^{-\Omega(t-s)}\|F(u(s))\|_{H^m}ds.
    \end{aligned}
  \end{equation}
  Using the bootstrap assumption
  \autoref{M-eq:section4-main-result:6}, namely
  $\|u(s)\|_{H^m}\leq \sqrt{2}E_m(0)$ for $s\in [0,t]$, yields
  \begin{equation}
    \label{eq:sec5:5b}
    \begin{aligned}
      \frac{\Omega^2}{\sqrt{2}}\int_0^t e^{-\Omega(t-s)}\|u(s)\|_{H^m}ds
      &\leq E_m(0)\Omega^2 \int_0^t e^{-\Omega(t-s)}ds \\
      &= E_m(0)\Omega \left( 1-e^{-\Omega t} \right).
    \end{aligned}
  \end{equation}
  Moreover, by the nonlinear bound \autoref{eq:sec5:14}, we conclude
  that
  \begin{equation}
    \label{eq:sec5:6b}
    \begin{aligned}
      \sqrt{2}\int_0^t e^{-\Omega(t-s)}\|F(u(s))\|_{H^m}ds
      &\leq \left( 1-e^{-\Omega t} \right) \frac{\sqrt{2}}{\Omega}\sup_{[0,t]} \|F(u(s))\|_{H^m} \\
      &\leq \left( 1-e^{-\Omega t} \right)\frac{\sqrt{2}}{\Omega}C(\delta)\sup_{[0,\infty)}\|a(s,\cdot)\|_{H^m}.
    \end{aligned}
  \end{equation}
  holds.
  Combining \autoref{eq:proof:9}–\autoref{eq:sec5:6b} gives
  \begin{equation}
    \label{eq:bug:1}
    \begin{aligned}
      E_m(t) \leq  e^{-\Omega t} E_m(0) +  \left( 1-e^{-\Omega t} \right)\left[E_m(0)\Omega + \tfrac{\sqrt{2}}{\Omega}C(\delta)\sup_{[0,\infty)}\|a(s,\cdot)\|_{H^m} \right].
    \end{aligned}
  \end{equation}
  We wish to show now that
  \begin{equation}
    \label{eq:section5-appendix:7}
    E_m(t)\leq \left(1 - \varepsilon'\right) E_m(0)
  \end{equation}
  holds for $t\in [T_1,T]$.
  By \autoref{eq:bug:1} it suffices to show that
  \begin{equation}
    \label{eq:bug:2}
    e^{-\Omega t} E_m(0) +  \left( 1-e^{-\Omega t} \right)\left[ 
      E_m(0)\Omega + \tfrac{\sqrt{2}}{\Omega}C(\delta)\sup_{[0,\infty)}\|a(s,\cdot)\|_{H^m} 
    \right] \leq \left(1 - \varepsilon'\right) E_m(0).
  \end{equation}
  is satisfied.
  Multiplying both sides of \autoref{eq:bug:2} by $e^{\Omega t}$ (recall
  $e^{\Omega t}<1$) we obtain
  \begin{equation}
    \label{eq:section5-appendix:8}
    \begin{aligned}
      E_m(0) + \left( e^{\Omega t}-1  \right)\left[ E_m(0)\Omega + \tfrac{\sqrt{2}}{\Omega}C(\delta)\sup_{[0,\infty)}\|a(s,\cdot)\|_{H^m} \right]
      \leq e^{\Omega t} \left(1 - \varepsilon'\right) E_m(0).
    \end{aligned}
  \end{equation}
  Rearranging terms in \autoref{eq:section5-appendix:8}, yields
  \begin{equation}
    \label{eq:section5-appendix:9}
    \begin{aligned}
      \left( e^{\Omega t}-1  \right)\tfrac{\sqrt{2}}{\Omega}C(\delta)\sup_{[0,\infty)}\|a(s,\cdot)\|_{H^m}
      \leq E_m(0) \left[e^{\Omega t}\left((1-\varepsilon')-\Omega\right) - (1-\Omega)\right].
    \end{aligned} 
  \end{equation}
  Note that by \ref{item:section5-appendix:1} the right hand side of
  \autoref{eq:section5-appendix:9} is strictly positive for $t\geq T_1$.
  Dividing by $(e^{\Omega t}-1)>0$ for $t\geq T_1$, we obtain
  \begin{equation}
    \label{eq:sec5:9} 
    \begin{aligned}
      &  \tfrac{\sqrt{2}}{\Omega}C(\delta)\sup_{[0,\infty)}\|a(s,\cdot)\|_{H^m} \\ 
      \leq  & E_m(0)\left[e^{\Omega t}\left(1-\varepsilon'-\Omega\right)-\left(1-\Omega 
           \right)\right]\left[e^{\Omega t}-1\right]^{-1}   = E_m(0)g(t).
    \end{aligned}
  \end{equation}
  Recall that $0<g(T_1)\leq g(t) $ for all $t\geq T_1$.
  Finally, choose $\varepsilon_{1}$ such that
  \begin{equation}
    \label{eq:section5-appendix:11}
    \sup_{[0,\infty)}\|a(s,\cdot)\|_{H^m} \leq \varepsilon_1<\tfrac{\Omega}{\sqrt{2}C(\delta)}g(T_1)E_m(0),
  \end{equation}
  is satisfied, which ensures the following inequality
  \begin{equation}
    \label{eq:section5-appendix:14}
    \tfrac{\sqrt{2}C(\delta)}{\Omega} \sup_{[0,\infty)}\|a(s,\cdot)\|_{H^m}
    \leq \tfrac{\sqrt{2}C(\delta)}{\Omega}\varepsilon_{1} < E_m(0)g(T_1)\leq E_m(0)g(t), \quad t\geq T_1.
  \end{equation}
  Thus \autoref{eq:sec5:9} holds and
  therefore~\autoref{eq:section5-appendix:7} is satisfied for $t\in [T_1,T]$. 
\end{proof}
\end{enumerate}

Having completed this step, the proof of \autoref{M-lem:section4-main-result:1} follows.

\renewenvironment{proof}[1][\proofname]{{\bfseries Proof of #1.}}{\hfill$\blacksquare$}
\end{proof}

\renewenvironment{proof}[1][\proofname]{{\bfseries Proof of #1.}}{\hfill$\blacksquare$}

\begin{proof}[~\fullref{M-lem:section4-main-result:2}]
  \label{sec:proof-fullr-lem:s-1}
  We begin decomposing $H^m =\mathbullet H^m\oplus \setR$, where
  $\mathbullet H^m$ denotes functions with zero mean and $\setR$ the constants. 
  This is an orthogonal decomposition under the inner product
  \begin{equation*}
    \langle f,g\rangle =\frac{1}{(2\pi)^3}\int_{\setT^3} f(x)\overline{g(x)}d^{3}x.
  \end{equation*}
  For any $u \in H^m$, we obtain the decomposition
  \begin{equation*}
    u=u_{h} + \bar u,
  \end{equation*}
  where the function $\bar{u}$ is given by
  \begin{equation*}
    \bar{u}=\overline{u(t)}=\langle 1, u\rangle=\frac{1}{(2\pi)^3}\int_{\setT^3} u(t,x) 
    d^{3}x,
  \end{equation*}
  and $u_h=u-\bar u\in \mathbullet H^m$. 
  It follows that
  \begin{equation*}
    \|u\|_{H^m}^2=\|u_h\|_{H^m}^2+|\bar u|^2
  \end{equation*}
  is true, since $u_{h}$ is orthogonal to $\bar u$.

  We compute each component of the norm, beginning with $\|u_h\|_{H^m}^2$. 
  Because $\partial_{\vec{\alpha}}u_{h}$ has zero mean for
  $|\vec{\alpha}|\leq m$, Wirtenger's inequality \cite{robinson2001} yields
  the following inequality
  \begin{equation}
    \label{eq:section5-appendix:1}
    \int\limits_{\setT^3}\left\vert  \partial_{\vec{\alpha}}u_{h} \right\vert^{2}d^3x
    \leq \int\limits_{\setT^3}\left\vert\nabla  \partial_{\vec{\alpha}}u_{h} \right\vert^{2}d^3x.
  \end{equation}
  Recall that $u\in H^{m+1}(\setT^{3})$ holds.
  Using \fullref{M-lem:section4-main-result:1} and the structure of
  the energy as defined by
  \autoref{M-eq:decaying-energy-estimates-nordstroem:3}, we conclude
  that
  \begin{equation}
    \label{eq:section5-appendix:2}
    \begin{aligned}
      \|u_{h}(T)\|_{H^m}^{2} 
      &  = \sum\limits_{|\alpha|\leq m}\left\Vert \partial_{\vec{\alpha}} u_{h}(T) 
        \right\Vert^{2}_{L^2} 
        \leq \sum\limits_{|\alpha|\leq m}\left\Vert \nabla \left(\partial_{\vec{\alpha}} 
        u_{h} \right)(T) \right\Vert^{2}_{L^2} \\ &
                                            =\sum\limits_{|\alpha|\leq m}\left\Vert \nabla \left( \partial_{\vec{\alpha}} 
                                            u\right)(T) \right\Vert^{2}_{L^2}\\ & \leq
                                                                           2 E_m^2(T)\leq2 \left(1-\varepsilon'\right)^2 E^2_{m}(0).
    \end{aligned}
  \end{equation}
  is true.
  It remains to estimate $ |\bar{u}(T)|^2$.  

  Since $u(t,x)$ satisfies the equation \makeatletter
  \def\plain@equationautorefname{equation}
  \begin{equation}
    \label{eq:strict-inquality:5}
    \partial_t^2 u- \Delta u+(2\Omega) \partial_t u_{t}=e^{-\kappa t} a(t, x)(1+u)^{\mu},
  \end{equation}
  we project \autoref{eq:strict-inquality:5} onto $\setR$ resulting in
  \begin{equation}
    \label{eq:section5-appendix:4}
    \langle \partial_t^2 u- \Delta u+(2\Omega) \partial_t u_{t},1 \rangle =  \partial_t^2 
    \bar u + (2\Omega) \partial_t \bar u_{t}+ \bar{F}(t),
  \end{equation}
  where $\bar F(t)$ is defined as:
  \begin{equation}
    \label{eq:section5-appendix:15}
    \bar F(t) = \left\langle   e^{-\kappa t} a(t, x)(1+u)^{\mu},1  \right\rangle=
    \left( \frac{1}{2\pi} \right)^{3}\int_{\setT^3}  e^{-\kappa t}a(t, x)(1+u)^{\mu}  d^3x.
  \end{equation}
  Recall $\bar u = \overline{u(t)}$ it satisfies an ordinary
  differential equation:
  \begin{equation}
    \label{eq:section5-appendix:13}
    \frac{d^2}{dt^2} \bar u+2\Omega \frac{d}{dt} \bar u=\bar F(t).
  \end{equation}

  By the assumptions of \fullref{M-thr:section4-main-result:1}, the
  initial data have zero mean, or in other words
  $\bar u(0)= \frac{d}{dt} \bar u (0)=0$.
  Hence the solution of \autoref{eq:section5-appendix:13} is given by
  \begin{equation}
    \label{eq:strict-inquality:9}
    \bar u(t)=\frac{1}{(2\Omega)} \int_{0}^{t}\left[\left(1-e^{-(2\Omega)(T-\tau)}\right) 
      \bar F(\tau)\right] d \tau, \qquad \tau\in [0,T].
  \end{equation}
  \makeatletter \def\plain@equationautorefname{Inequality}

  To estimate $|\bar F(\tau)|$, note that by the bootstrap assumption
  $\|u(\tau)\|^{2}_{H^m}\leq 2E_m^{2}(0)$ holds for $\tau\in [0,T]$ and by
  \autoref{eq:sec5:13} and
  \autoref{M-eq:global-existence-lavi-bootstrap-new:70}, we obtain
  \begin{equation}
    \label{eq:section5-appendix:19}
    \begin{aligned}
      \big|\bar{F}(\tau)\big|
      & =\frac{1}{(2\pi)^3}\bigg|\int_{\setT^3}e^{-\kappa 
        \tau}a(t, x)(1+u(\tau))^{\mu} d^{3}x\bigg| \\
      & \leq C(\delta) e^{-\kappa \tau}\sup_{[0,\infty)}\|a(t, \cdot)\|_{H^m}\leq  C(\delta) e^{-\kappa \tau}\varepsilon_{2},
    \end{aligned}
  \end{equation}
  provided that
  \begin{equation}
    \sup_{[0,\infty)}\|a(t, \cdot)\|_{H^m}\leq \varepsilon_{2}
  \end{equation}
  holds.

  Therefore we can estimate $ \big| \bar u(T)\big|$ as follows:
  \begin{equation}
    \label{eq:section5-appendix:12}
    \begin{aligned}
      \big| \bar u(T)\big|
      & \leq \frac{\varepsilon_2C(\delta)}{2\Omega} 
        \int_{0}^{T}\left[\left(1-e^{-(2\Omega)(T-\tau)}\right) e^{-\kappa \tau}\right] d \tau
        \leq \frac{\varepsilon_2C(\delta)}{2\Omega}  \int_{0}^{T}\left[e^{-\kappa \tau}\right] d \tau \\
      & \leq \frac{\varepsilon_2C(\delta)}{2\Omega \kappa}\left( 1-e^{-\kappa T} \right) 
        \leq \frac{\varepsilon_2C(\delta)}{2\Omega \kappa}.
    \end{aligned}
  \end{equation}
  Choosing $\varepsilon_{2}$ so that
  \begin{equation}
    \label{eq:section5-appendix:18}
    \left( \frac{\varepsilon_2C(\delta)}{2\Omega \kappa}  \right)^{2}\leq 2 \left( 
      \varepsilon^{\prime}-\tfrac{3}{4} \left( \varepsilon^{\prime} \right)^{2} \right)E_m^2(0),
  \end{equation}
  holds and also $\varepsilon_2\leq \varepsilon_1$ is satisfied, we then combine it with
  \autoref{eq:section5-appendix:2}, and obtain the following
  inequality
  \begin{equation}
    \label{eq:section5-appendix:10}
    \begin{aligned}
      \|u(T)\|_{H^m}^2 
      &= \|u_{h}(T)\|_{H^m}^2+ |\bar u|^{2} \\
      &\leq 2 \left(1-\varepsilon'\right) E^2_{m}(0)+ 2 \left( \varepsilon^{\prime}-\tfrac{3}{4} \left( \varepsilon^{\prime} \right)^{2} \right)E_m^2(0) \\
      &=2E_m^2(0)\left( 1-\tfrac{\varepsilon^{\prime}}{2} \right)^2 < 2E_m^2(0).
    \end{aligned}
  \end{equation}
  which is the desired inequality.
\end{proof}

\begin{rem}[On the choice of $\varepsilon_{2}$]
  \label{rem:section5-appendix:1}
  The choice of $\varepsilon_2$ depends on $\varepsilon^{\prime}$ and on
  $\varepsilon_1 $, and therefore it depends on the initial data.
\end{rem}

\renewenvironment{proof}[1][\proofname]{{\bfseries Proof of #1.}}{\hfill$\blacksquare$}

\section{Asymptotic behaviour}
\label{sec:asymptotic-behaviour}
  The main result, \fullref{M-thr:section4-main-result:1}, ensures the
  existence of a global solution to the Cauchy problem
  \eqref{M-eq:wave:1}--\eqref{M-eq:wave:2} that remains small for all
  $t \geq 0$. In this section we show that the solution converges to a
  constant, denoted by $c_0$.
\makeatletter \def\plain@equationautorefname{Equation}
\makeatother

At first glance, one might expect the global solution to decay to zero
as $t \to \infty$. 
Yet a closer look at \autoref{M-eq:section5-appendix:13}, together
with the estimate for $\bar F(t)$ in \makeatletter
\def\plain@equationautorefname{Inequality} \makeatother
\autoref{M-eq:section5-appendix:19}, suggests otherwise. 
Indeed, from the explicit form of the solution given in \makeatletter
\def\plain@equationautorefname{Equation} \makeatother
\autoref{M-eq:strict-inquality:9}, it follows that $\bar u(t)$ cannot
converge to zero, even though $\lim_{t\to\infty} \bar F(t)=0$. 
The strongest conclusion available, therefore, is that the solution
converges instead to a constant.

\begin{thm}[Asymptotic behaviour of the global solution]
  \label{thr:section6-asymptotics:1}
  Let $u$ be the global solution to the Cauchy problem
  \eqref{M-eq:wave:1}--\eqref{M-eq:wave:2}. 
  Then
  \begin{equation}
    \lim_{t\to \infty} u(t,x)= c_0
  \end{equation}
  holds, where \(c_0\) is defined by
  \begin{equation}
    c_0 = \lim_{t\to \infty} \left(\frac{1}{2\pi}\right)^3 \int_{\setT^3} u(t,x)\, d^{3}x.
  \end{equation}
\end{thm}

\begin{rem}[On the asymptotic behaviour]
 \label{rem:section6-asymptotics:2}
Recall that in Nordström’s theory of gravitation the metric $g_{\alpha\beta}$ is conformal
to the Minkowski metric $\eta_{\alpha \beta}$. Explicitly written as
\begin{equation}
  \label{eq:section6-asymptotics:2}
  g_{\alpha\beta}=  \phi^2(t,x)\,\eta_{\alpha \beta}.
\end{equation}
The background metric in our case is $e^{2\Omega t}\eta_{\alpha \beta}$. 
Thus for large $t$, the asymptotic form of the metric is
\begin{equation}
  \label{eq:section6-asymptotics:3}
  g_{\alpha\beta}=  e^{2\Omega t}(1+c_0)^2\eta_{\alpha \beta},
\end{equation}
which is a slight perturbation of the background geometry.
\end{rem}

  \begin{rem}[About the choice of the Energy functional]
    We observe that the energy used to prove global existence, namely
    \autoref{M-eq:decaying-energy-estimates-nordstroem:3}, is unsuited
    for establishing the asymptotic behaviour. 
    Instead, we adopt the standard energy of the wave equation:
    \begin{equation}
      \label{eq:section6-asymptotics:13}
      \begin{aligned}
        \mathds{E}^2(t) &= \frac{1}{2} \sum_{0 \leq |\vec{\alpha}| \leq m} 
                          \int_{\setT^3}\left(|\partial_t \partial_{\vec{\alpha}} u(t,x)|^2 + |\nabla ( 
                          \partial_{\vec{\alpha}} u(t,x))|^2\right) d^{3}x  \\
                        &= \frac{1}{2}\left(\|\partial_t u(t)\|_{H^m}^2 + \|\nabla 
                          u(t)\|_{H^m}^2\right).
      \end{aligned}
    \end{equation}
\end{rem}

\subsection{Proof of \autoref{thr:section6-asymptotics:1}}
\label{sec:proof-autor-asympt}
\hfill

\textsc{Proof Sketch:} We will use the energy
\autoref{eq:section6-asymptotics:13}, together with the properties of
the source term in \autoref{M-eq:sec5:12}, to show that
\begin{equation}
  \label{eq:section6-asymptotics:1}
  \lim_{t\to\infty}\|\partial_t u(t)\|_{H^m}=0
\end{equation}
holds. 
The limit function then depends only on spatial variables, and must be
harmonic on the torus $\setT^3$. 
Such a function is necessarily constant.

\begin{proof}[~\autoref{thr:section6-asymptotics:1}]
  Since $u$ is the global solution, it satisfies
  \begin{equation}
    \label{eq:sec6:9}
    \partial_t^2 u +2\Omega \partial_t u -\Delta u = F(u(t)), \qquad t\geq 0, 
  \end{equation}
  where $F(u(t))$ is defined by \makeatletter
  \def\plain@equationautorefname{Equation} \makeatother
  \autoref{M-eq:sec5:12}. 
  Differentiating $\mathds{E}^2(t)$, given by
  \autoref{eq:section6-asymptotics:13}, with respect to $t$ yields
  \begin{equation}
    \label{eq:section6-asymptotics:14}
    \begin{aligned}
      \frac{d\mathds{E}^2}{dt}
      &= \sum_{0 \leq |\vec{\alpha}| \leq m}\int_{\setT^3}\partial_t\partial_{\vec{\alpha}} u\Big[\partial_t^2\partial_{\vec{\alpha}} u- \Delta(\partial_{\vec{\alpha}} u)\Big]\,d^{3}x \\
      &= -2\Omega \sum_{0 \leq |\vec{\alpha}| \leq m}\int_{\setT^3}\big(\partial_t\partial_{\vec{\alpha}} u\big)^2\,d^{3}x
        + \sum_{0 \leq |\vec{\alpha}| \leq m}\int_{\setT^3}\big(\partial_t\partial_{\vec{\alpha}} u\big)\,\partial_{\vec{\alpha}}(F(u))\,d^{3}x \\
      &= -2\Omega \sum_{0 \leq |\vec{\alpha}| \leq m}\int_{\setT^3}\Big[\big(\partial_t\partial_{\vec{\alpha}} u\big)^2
        + |\nabla(\partial_{\vec{\alpha}} u)|^2\Big]\,d^{3}x + 2\Omega \sum_{0 \leq |\vec{\alpha}| \leq m}\int_{\setT^3} |\nabla(\partial_{\vec{\alpha}} u)|^2\, d^{3}x \\
      &\quad + \sum_{0 \leq |\vec{\alpha}| \leq m}\int_{\setT^3}\big(\partial_t\partial_{\vec{\alpha}} u\big)\,\partial_{\vec{\alpha}}(F(u))\,d^{3}x \\
      & \leq -4\Omega \mathds{E}^2 + 2\Omega \|\nabla u\|_{H^m}^2 + \|\partial_t u\|_{H^m}\|F(u)\|_{H^m}.
    \end{aligned}
  \end{equation}
  The last step uses the Cauchy–Schwarz inequality. 
  From the definition of $F(u(t))$ in \makeatletter
  \def\plain@equationautorefname{Equation}
  \makeatother\autoref{M-eq:sec5:12}, together with \makeatletter
  \def\plain@equationautorefname{Inequality}
  \makeatother\autoref{M-eq:sec5:14}, it follows that
  \begin{equation}
    \label{eq:sec6:7}
    \|F(u(t))\|_{H^m}\leq C_1 e^{-\kappa t}, \qquad t\geq 0
  \end{equation}
  holds.
  Moreover, by \fullref{M-def:decaying-energy-estimates-nordstroem:2},
  and the fact that  this energy functional is abounded for all $ t$, expressed by
  \autoref{M-eq:global-existence-lavi-bootstrap-new:74} we conclude
  the following inequality
  \begin{equation}
    \label{eq:section6-asymptotics:6}
    \|\partial_t u(t)\|_{H^m} \leq 4E(t) \leq C_2, \qquad t\geq 0,
  \end{equation}
  is satisfied, which implies the inequality
  \begin{equation}
    \label{eq:section6-asymptotics:5}
    \frac{d\mathds{E}^2(t)}{dt} \leq -4\Omega \mathds{E}^2(t) + 2\Omega \|\nabla u(t)\|_{H^m}^2 + C_1C_2 e^{-\kappa t}.
  \end{equation}
  Applying the standard Gronwall inequality \cite[Appendix B]{evans99:_partial_differ_equat} to
  \autoref{M-eq:section6-asymptotics:5} results in the following
  inequality
  \begin{equation}
    \label{eq:section6-asymptotics:4}
    \mathds{E}^2(t)\leq e^{-4\Omega(t-t_0)}\mathds{E}^2(t_0)+2\Omega \int_{t_0}^t 
    e^{-4\Omega(t-\tau)}\|\nabla u(\tau)\|_{H^m}^2 d\tau + C_1 C_2 \int_{t_0}^t 
    e^{-4\Omega(t-\tau)}e^{-\kappa \tau} d\tau  
  \end{equation}
  that holds for all $t\geq t_0$.

  Since the solution $u$ exists globally and remains bounded, we
  define
  \begin{equation}
    \label{eq:section6-asymptotics:7}
    \mu^2 = \limsup_{t\to \infty} \|\nabla u(t)\|_{H^m}^2.
  \end{equation}
  For any $\epsilon >0$, there exists $t_0>0$ such that
  \begin{equation}
    \label{eq:section6-asymptotics:8}
    \|\nabla u(t)\|_{H^m}^2\leq \mu^2 +\epsilon
  \end{equation}
  holds for all $t\geq t_0$. 
  Substituting this into \autoref{eq:section6-asymptotics:4}, we
  obtain
  \begin{equation}
    \label{eq:sec6:2}
    \mathds{E}^2(t)\leq e^{-4\Omega(t-t_0)}\mathds{E}^2(t_0)+2\Omega \int_{t_0}^t 
    e^{-4\Omega(t-\tau)}(\mu^2 +\epsilon) d\tau + C_1 C_2 \int_{t_0}^t 
    e^{-4\Omega(t-\tau)}e^{-\kappa \tau} d\tau.
  \end{equation}

  We evaluate the two integrals in \autoref{eq:sec6:2}. 
  We estimate the first integral as follows:
  \begin{equation}
    \label{eq:sec6:3}
    \begin{aligned}
      2\Omega \int_{t_0}^t e^{-4\Omega(t-\tau)}(\mu^2 +\epsilon) d\tau 
      &=\tfrac{1}{2}(\mu^2 +\epsilon) e^{-4\Omega t}\big(e^{4\Omega t}-e^{4\Omega t_0}\big) \\
      &=\tfrac{1}{2}(\mu^2 +\epsilon)\big(1-e^{-4\Omega (t-t_0)}\big) \\
      &\leq \tfrac{1}{2}(\mu^2 +\epsilon).
    \end{aligned}
  \end{equation}
  For the second integral we use the inequality $4\Omega -\kappa>0$,
  which allows us to estimate that integral as follows.
  \begin{equation}
    \label{eq:sec6:4}
    \begin{aligned}
      \int_{t_0}^t  e^{-4\Omega t} e^{(4\Omega-\kappa)\tau} d\tau
      &= e^{-4\Omega t}\frac{1}{4\Omega -\kappa}\Big(e^{(4\Omega -\kappa) t}-e^{(4\Omega -\kappa)t_0}\Big) \\
      &= \frac{1}{4\Omega -\kappa }\Big(e^{-\kappa t}-e^{-4\Omega(t-t_0)}e^{-\kappa t_0}\Big) \\
      &\leq \frac{e^{-\kappa t}}{4\Omega -\kappa }.
    \end{aligned}
  \end{equation}
  If $4\Omega -\kappa\leq 0$, then
  \begin{equation}
    \label{eq:sec6:4b}
    \int_{t_0}^t  e^{-4\Omega t} e^{(4\Omega-\kappa)\tau} d\tau 
    \leq e^{-4\Omega t} \int_{t_0}^t d\tau 
    \leq e^{-4\Omega t}(t-t_0).
  \end{equation}
  is true.
  In either cases, we can conclude that
  \begin{equation}
    \label{eq:section6-asymptotics:10}
    \limsup_{t\to \infty} \int_{t_0}^t e^{-4\Omega(t-\tau)}e^{-\kappa \tau} d\tau =0
  \end{equation}
  holds.
  Thus, by \autoref{eq:section6-asymptotics:4} and
  \autoref{eq:sec6:2}, we conclude that
  \begin{equation}
    \label{eq:section6-asymptotics:9}
    \lim_{t\to\infty}  \mathds{E}^2(t)\leq  \tfrac{1}{2}(\mu^2 +\epsilon)
  \end{equation}
  is satisfied.
  It then follows that
  \begin{equation}
    \label{eq:section6-asymptotics:15}
    \begin{aligned}
      \limsup_{t\to\infty}  \mathds{E}^2(t)
      &=\limsup_{t\to\infty}\tfrac{1}{2}\Big(
        \|\partial_t u(t)\|_{H^m}^2+\|\nabla u(t)\|_{H^m}^2\Big) \\
      &=\limsup_{t\to\infty}\tfrac{1}{2}
        \|\partial_t u(t)\|_{H^m}^2 +\tfrac{1}{2}\mu^2 \\
      &\leq \tfrac{1}{2}(\mu^2+\epsilon).
    \end{aligned}
  \end{equation}
  is true.
  Therefore we conclude
  \begin{equation}
    \label{eq:section6-asymptotics:16}
    \limsup_{t\to\infty} \|\partial_t u(t)\|_{H^m}^2 \leq \epsilon
  \end{equation}
  holds for any $\epsilon>0$. 
  That is why we conclude that
  \begin{equation}
    \label{eq:section6-asymptotics:17}
    \limsup_{t\to\infty} \|\partial_t u(t)\|_{H^m} =0,
  \end{equation}
  is satisfied  and then by \fullref{M-prop:Sobolev-inequalyi}, that
  \begin{equation}
    \label{eq:sec6:8}
    \limsup_{t\to\infty} \|\partial_t u(t)\|_{C^1} =0 \qquad     \limsup_{t\to\infty} \|\partial_{tt} u(t)\|_{C^0} =0 \qquad 
  \end{equation}
 is true.

 Since the limit function does not depend on the time $ t $, we set
  \begin{equation*}
    \widetilde u(x):= \lim_{t\to \infty}u(t,x).
  \end{equation*}
  Then, using \autoref{eq:sec6:7}, \makeatletter
  \def\plain@equationautorefname{Equation} \makeatother (together,
  with the Sobolev inequality), \autoref{eq:sec6:8}, and
  \autoref{eq:sec6:9} we observe that $\widetilde u$ satisfies the
  equation
  \begin{equation}
    -\Delta \widetilde u  =0.
  \end{equation}
  Its Fourier expansion,
  \begin{equation*}
    \widetilde u(x)=\sum_{\vec{n}\in \setZ^3} c_{\vec{n}}e^{i(\vec{n}\cdot x)},
  \end{equation*}
  yields
  \begin{equation*}
    \Delta \widetilde u(x)=-\sum_{\vec{n}\in \setZ^3} c_{\vec{n}}|\vec{n}|^2e^{i(\vec{n}\cdot x)}=0,
  \end{equation*}
  implying that $c_{\vec{n}}=0$ is satisfied for all $\vec{n}\neq 0$. 
  Thus we conclude that
  \begin{equation*}
    \widetilde u(x) =c_0=\left(\tfrac{1}{2\pi}\right)^3\int_{\setT^3}  \widetilde u(x)\,d^3x
    =\lim_{t\to \infty}\left(\tfrac{1}{2\pi}\right)^3\int_{\setT^3}  u(t,x)\,d^3x
  \end{equation*}
  is true.

\end{proof}

\printbibliography
\end{document}